\documentclass[12pt,reqno]{amsart}

\usepackage{amscd}
\usepackage{epsfig}
\usepackage{float}
\usepackage{url}
\usepackage[all]{xy}
\newcommand{\N}{{\mathbb N}}
\newcommand{\Z}{{\mathbb Z}}
\newcommand{\Q}{{\mathbb Q}}
\newcommand{\C}{{\mathbb C}}

\newcommand{\GG}{{\mathcal G}}

\newcommand{\www}{\widetilde}

\newcommand{\oooo}{\overline}

\newcommand{\paa}{\partial}

\DeclareMathOperator{\Aut}{Aut}

\DeclareMathOperator{\divis}{div}

\DeclareMathOperator{\id}{id}

\DeclareMathOperator{\lcm}{lcm}

\DeclareMathOperator{\ord}{ord}

\DeclareMathOperator{\supp}{supp}

\DeclareMathOperator{\tr}{tr}

\begin{document}

\theoremstyle{plain}
\newtheorem{lemma}{Lemma}[section]
\newtheorem{definition/lemma}[lemma]{Definition/Lemma}
\newtheorem{theorem}[lemma]{Theorem}
\newtheorem{proposition}[lemma]{Proposition}
\newtheorem{corollary}[lemma]{Corollary}
\newtheorem{conjecture}[lemma]{Conjecture}
\newtheorem{conjectures}[lemma]{Conjectures}

\theoremstyle{definition}
\newtheorem{definition}[lemma]{Definition}
\newtheorem{withouttitle}[lemma]{}
\newtheorem{remark}[lemma]{Remark}
\newtheorem{remarks}[lemma]{Remarks}
\newtheorem{example}[lemma]{Example}
\newtheorem{examples}[lemma]{Examples}

\title[Seven combinatorial problems around qh singularities]
{Seven combinatorial problems around quasihomogeneous singularities} 

\author{Claus Hertling and Philip Zilke}

\address{Claus Hertling\\
Lehrstuhl f\"ur Mathematik VI, Universit\"at Mannheim, Seminargeb\"aude
A 5, 6, 68131 Mannheim, Germany}

\email{hertling@math.uni-mannheim.de}

\address{Philip Zilke\\
Lehrstuhl f\"ur Mathematik VI, Universit\"at Mannheim, Seminargeb\"aude
A 5, 6, 68131 Mannheim, Germany}

\email{Philip.Zilke@gmx.de}

\date{January 25, 2018}

\subjclass[2010]{32S40, 12Y05, 05C22, 05C25}

\keywords{quasihomogeneous singularity, weight system, monodromy, 
characteristic polynomial, combinatorial problems, Orlik blocks}

\thanks{This work was supported by the DFG grant He2287/4-1
(SISYPH)}


\begin{abstract}
{\Small 
This paper proposes seven combinatorial problems around formulas for the 
characteristic polynomial and the spectral numbers of a quasihomogeneous 
singularity. One of them is a new conjecture on the characteristic 
polynomial. It is an amendment to an old conjecture of Orlik on the
integral monodromy of a quasihomogeneous singularity.
The search for a combinatorial proof of the new conjecture led us to 
the seven purely combinatorial problems.}
\end{abstract}

\maketitle

\tableofcontents

\setcounter{section}{0}

\section{Introduction}\label{s1}
\setcounter{equation}{0}

\noindent
This paper proposes seven combinatorial problems around formulas for the 
characteristic polynomial and the spectral numbers of a quasihomogeneous 
singularity. One of them is a new conjecture on the characteristic 
polynomial. It is an amendment to an old conjecture of Orlik on the
integral monodromy of a quasihomogeneous singularity.
The search for a combinatorial proof of the new conjecture led us to 
the seven purely combinatorial problems.

We start with a result on $\Z$-lattices with automorphisms.
Then we describe Orlik's conjecture and our new conjecture.
Finally, we give a rough outline of the seven problems. 

\begin{definition}\label{t1.1}
Let $M\subset\N=\{1,2,3,...\}$ be a finite nonempty subset. 
Its {\it Orlik block} is a pair $(H_M,h_M)$ with 
$H_M$ a $\Z$-lattice of rank $\sum_{m\in M}\varphi(m)$ and 
$h_M:H_M\to H_M$ an automorphism with characteristic polynomial
$\prod_{m\in M}\Phi_m$ ($\Phi_m$ is the $m$-th cyclotomic polynomial)
and with a cyclic generator $e_1\in M$, i.e.
\begin{eqnarray}\label{1.1}
H_M =\bigoplus_{j=1}^{rk M}\Z\cdot h_M^{j-1}(e_1).
\end{eqnarray}
$(H_M,h_M)$ is unique up to isomorphism.
$\Aut_{S^1}(H_M,h_M)$ denotes the group of all automorphisms of $H_M$
which commute with $h_M$ and which have all eigenvalues in $S^1$.
\end{definition}

Definition \ref{t6.1} below enriches the set $M$ to a directed graph $\GG(M)$.
An edge goes from $m_1\in M$ to $m_2\in M$ if $\frac{m_1}{m_2}$ is a
power of a prime number $p$. Then it is called a $p$-edge.
The main result in \cite{He2} is cited precisely in 
theorem \ref{t6.2} below. Roughly it is as follows.

\begin{theorem}\label{t1.2} \cite[Theorem 1.2]{He2}
Let $(H_M,h_M)$ be the Orlik block of a finite nonempty subset $M\subset\N$.
Then $\Aut_{S^1}(H_M,h_M)=\{\pm h_M^k\, |\, k\in\Z\}$ if and only if
condition (I) or condition (II) in theorem \ref{t6.2}
are satisfied. They are conditions on the graph $\GG(M)$.
\end{theorem}

A weight system ${\bf w}=(w_1,...,w_n)$ with $w_i\in\Q_{>0}$
equips any monomial ${\bf x}^{\bf j}=x_1^{j_1}...x_n^{j_n}$ with a 
weighted degree $\deg_{\bf w}{\bf x}^{\bf j}:=\sum_{i=1}^n w_ij_i$. 
A polynomial $f\in\C[x_1,...,x_n]$ is a {\it quasihomogeneous singularity}
if for some weight system ${\bf w}$ with $w_i\in\Q\cap (0,1)$
each monomial in $f$ has weighted degree 1 and if the functions 
$\frac{\paa f}{\paa x_1},...,\frac{\paa f}{\paa x_n}$ vanish simultaneously
only at $0\in\C^n$. Then the {\it Milnor lattice} 
$H_{Milnor}:=H_{n-1}(f^{-1}(1),\Z)$ is a 
$\Z$-lattice of some rank $\mu\in\N$ \cite{Mi}, which is called
{\it Milnor number}. It comes equipped with a natural automorphism
$h_{mon}:H_{Milnor}\to H_{Milnor}$ of finite order, the {\it monodromy}.
Thus its characteristic polynomial has the form 
\begin{eqnarray*}
p_{ch,h_{mon}} &= \prod_{m\in M_1}\Phi_m^{\nu(m)}
\end{eqnarray*}
for a finite subset $M_1\subset \N$ and a function $\nu:M_1\to\N$.
Denote $\nu_{max}:=\max(\nu(m)\, |\, m\in\N)$ and for $j=1,...,s_{max}$
\begin{eqnarray*}
M_j:=\{m\in M_1\, |\, \nu(m)\geq j\},\quad
g_j:=\prod_{m\in M_j} \Phi_m.
\end{eqnarray*}
Then 
\begin{eqnarray*}
M_1\supset M_2\supset ...\supset M_{\nu_{max}}\neq\emptyset\\
\textup{and}\qquad p_{ch,h_{mon}} =\prod_{j=1}^{\nu_{max}} g_j.
\end{eqnarray*}
The polynomials $g_1,...,g_{\nu_{max}}$ are called {\it elementary 
divisors} of $p_{ch,h_{mon}}$. 

\begin{conjecture}\label{t1.3}
(Orlik's conjecture, \cite[conjecture 3.1]{Or}) 
For any quasihomogeneous singularity, there is an isomorphism
\begin{eqnarray*}
(H_{Milnor},h_{mon})\cong \bigoplus_{j=1}^{\nu_{max}} (H_{M_j},h_{M_j}).
\end{eqnarray*}
\end{conjecture}

The conjecture is known to be true for curve singularities \cite{MW}
and a few other cases, but it is still (after 45 years)
open in general. The following conjecture is an amendment to Orlik's
conjecture. 

\begin{conjecture}\label{t1.4}
For any quasihomogeneous singularity, each of the sets $M_1,...,M_{\nu_{max}}$
satisfies condition (I) in theorem \ref{t6.2}.
\end{conjecture}

If this is true for some singularity, 
then theorem \ref{t1.2} gives for these sets
$\Aut_{S^1}(H_{M_j},h_{M_j})\cong\{\pm h_{M_j}^k\, |\, k\in\Z\}$.
If also Orlik's conjecture holds, then 
this is helpful in determining the automorphisms of the Milnor lattice
which respect the monodromy (and intersection form or Seifert form).

Orlik's conjecture \ref{t1.3} concerns the Milnor lattice. 
Any proof requires to go into geometry. 
But we hope that our conjecture \ref{t1.4} is amenable to a combinatorial proof.
It just concerns the characteristic polynomial.
Milnor and Orlik \cite{MO} proved a formula which expresses this  
in terms of the weight system ${\bf w}$ of the quasihomogeneous
singularity. It says $\divis p_{ch,h_{mon}}=D_{\bf w}$, where
$D_{\bf w}$ is defined in \eqref{3.9}. See theorem \ref{t3.9}.
Therefore we hope that there will be a purely combinatorial
proof of conjecture \ref{t1.4} dealing solely with properties of
${\bf w}$. This is problem 6 below. For most of the other problems,
we need two more data.

First, a quasihomogeneous singularity comes also equipped with 
{\it exponents} $\alpha_1,...,\alpha_\mu\in \Q\cap (0,n)$.
They are slightly finer invariants than $p_{ch,h_{mon}}$. They satisfy
\begin{eqnarray*}
\divis p_{ch,h_{mon}} =\sum_{j=1}^\mu \langle e^{2\pi i\alpha_j}\rangle 
\end{eqnarray*}
and, for any $d\in \N$ with $v_i:=d\cdot w_i\in \N$ for all $i\in\{1,...,n\}$,
\begin{eqnarray*}
\sum_{j=1}^\mu t^{d\cdot\alpha_j} = \rho_{{\bf v},d}
\end{eqnarray*}
where $\rho_{{\bf v},d}$ is defined in \eqref{3.8}.
See theorem \ref{t3.9}.

Second, the weight systems ${\bf w}$ for which quasihomogeneous singularities
exist, can be characterized by a combinatorial condition $(C1)$
(and equivalent combinatorial conditions $(C1)'$ and $(C2)$, see lemma
\ref{t3.3}). This is cited in theorem \ref{t3.5}. 
It was proved first by Kouchnirenko \cite[Remarque 1.13 (i)]{Ko1}. 
The necessity of $(C1)$ had
already been seen by K. Saito \cite{Sa1}, the sufficiency not.
A weaker combinatorial property $\oooo{(C1)}$ is equivalent to 
$\rho_{{\bf v},d}\in\Z[t]$ \eqref{3.16}. 

The seven problems are given in detail in the later chapters. Roughly,
they are as follows.

\begin{list}{}{}
\item[{\bf Problem 1:}] (Remark \ref{t3.8} (i))
Let $({\bf v},d)=(v_1,...,v_n,d)\in\N^{n+1}$ 
with $d>\max_i v_i$ be given which satisfies $\oooo{(C1)}$.
Write $\rho_{({\bf v},d)}=\sum_{\alpha\in\frac{1}{d}\Z}\sigma(\alpha)\cdot
t^{d\cdot\alpha}\in\Z[t]$. Is 
$D_{\bf w}=\sum_{\alpha}\sigma(\alpha)\cdot\langle e^{2\pi i\alpha}\rangle$?
\item[{\bf Problem 2 :}] (Remark \ref{t3.11} (ii))
Let $({\bf v},d)=(v_1,...,v_n,d)\in\N^{n+1}$ 
with $d>\max_i v_i$ be given which satisfies $(C1)$. 
Give combinatorial proofs of the formulas in theorem \ref{t3.9} 
which connect $D_{\bf w}$ and $\rho_{({\bf v},d)}$ with the exponents
and with one another.
\item[{\bf Problem 3:}] (Remark \ref{t3.11} (iii))
Make some good use of the conditions for $J$ with $|J|\geq 2$ in $(C1)$.
\item[{\bf Problem 4:}] (Remarks \ref{t5.2}) 
Find examples different from Ivlev's example for weight systems 
${\bf w}$ which satisfy $\oooo{(C1)}$, but not $(C1)$. 
\item[{\bf Problem 5:}] (Remark \ref{t5.6})
Prove or disprove K. Saito's conjecture \ref{t5.4}
that $d_{\bf w}\in M_1$ or $\frac{d_{\bf w}}{2}\in M_1$ for ${\bf w}$ with $(C1)$.
Here $d_{\bf w}:=\lcm(\textup{denominator of }w_i\, |\, i\in\{1,...,n\})$.
\item[{\bf Problem 6:}] (Remark \ref{t6.4} (i))
Prove (or disprove) combinatorially conjecture \ref{t1.4}.
\item[{\bf Problem 7:}] (Remark \ref{t6.4} (iii))
Find a natural condition on products $f$ of cyclotomic polynomials which
implies for any elementary divisor of $f$ condition (I) in theorem \ref{t6.2}
and which is stable under tensor product. Prove that $D_{\bf w}$ for 
${\bf w}$ with $(C1)$ satisfies it
(this would prove conjecture \ref{t1.4}).
\end{list} 

Some comments: The problems 1, 2, 3 and 7 are motivated by problem 6, 
i.e. the wish to prove {\it combinatorially} conjecture \ref{t1.4}.
The problems 1 and 2 are closely related. 
A positive solution to one of them will probably also give 
a positive solution to the other one.
\cite{HK1} made good use of the conditions for $|J|=1$ in $(C1)$.
They give rise to a graph. But the problems here probably require
to involve also the conditions for $|J|\geq 2$. Problem 3 is vague,
but fundamental.
It looks surprisingly difficult to find solutions for the 
very concrete problem 4.
Problem 5 is motivated by the (more important) problems 6 and 7.
They are closely related. A positive solution of problem 6
goes probably via a positive solution of problem 7.

The paper is structured as follows.
Section \ref{s2} gives notations and basic facts 
around cyclotomic polynomials.
Section \ref{s3} introduces for abstract weight systems
${\bf w}$ and $({\bf v},d)$, the objects $D_{\bf w}$ and $\rho_{({\bf v},d)}$,
the conditions $\oooo{(C1)}$ and $(C1)$, and it states elementary facts
as well as the formulas and facts which hold for the weight systems ${\bf w}$
of quasihomogeneous singularities. This is all classical.
Section \ref{s4} gives more explicit formulas in the cases of the quasihomogeneous
singularities of cycle type and chain type. This builds on section \ref{s3} 
and on \cite{HK1} and is elementary.
Section \ref{s5} presents examples. Especially, it gives counterexamples
to the part of K. Saito's conjecture \ref{t5.4} which says that $d_{\bf w}\in M_1$ 
in the case of a weight system ${\bf w}$ with all $w_i<\frac{1}{2}$.
These counterexamples are interesting also in section \ref{s6}.
Section \ref{s6} formulates in theorem \ref{t6.2}
the main result from \cite{He2} on automorphisms of Orlik blocks.
It discusses conjecture \ref{t1.4}, it gives examples, 
and it proves conjecture \ref{t1.4} in special cases, which include
the cycle type, the chain type, the cases with $n=2$ and many of the
cases with $n=3$ (theorem \ref{t6.9}).

\section{Notations around cyclotomic polynomials}\label{s2}
\setcounter{equation}{0}

\noindent
This section fixes some notations and recalls some well known
formulas around products of cyclotomic polynomials.

In this paper $\N=\{1,2,3,...\}$ and $\N_0=\{0,1,2,3,...\}$.
Whenever a number $n\in\N$ is fixed then 
$N:=\{1,...,n\}$.

Denote by $S^{UR}\subset S^1$ the set of all unit roots.
Denote by $\Q\langle S^{UR}\rangle$ and $\Z\langle S^{UR}\rangle$
the group rings with elements 
$\sum_{j=1}^l b_j\langle\zeta_j\rangle$
where $b_j\in\Q$ respectively $b_j\in \Z$ 
and where $\zeta_j\in S^{UR}$,
with multiplication $\langle\zeta_1\rangle\cdot
\langle\zeta_2\rangle=\langle \zeta_1\cdot\zeta_2\rangle$.
The unit element is $\langle 1\rangle$.
The trace of an element $\sum_{j=1}^l b_j\langle\zeta_j\rangle$
is
\begin{eqnarray}\label{2.1}
\tr\left(\sum_{j=1}^l b_j\langle\zeta_j\rangle\right):=
\sum_{j=1}^l b_j\cdot\zeta_j\in\C.
\end{eqnarray}
The degree of it is 
\begin{eqnarray}\label{2.2}
\deg\left(\sum_{j=1}^l b_j\langle\zeta_j\rangle\right):=
\sum_{j=1}^l b_j\in\Q.
\end{eqnarray}
The trace map $\tr:\Q\langle S^{UR}\rangle\to\C$ 
and the degree map $\deg:\Q\langle S^{UR}\rangle\to\Q$
are ring homomorphisms.

The divisor of a unitary 
polynomial $f=(t-\lambda_1)\cdot ...\cdot (t-\lambda_l)\in\C[t]$
with $\lambda_j\in S^{UR}$ is 
\begin{eqnarray}\label{2.3}
\divis f:=\langle \lambda_1\rangle+...+\langle\lambda_l\rangle.
\end{eqnarray}
Of course $\tr(\divis f)=\lambda_1+...+\lambda_l$
and $\deg(\divis f)=\deg f$.

For two polynomials, $f$ as above and 
$g=(t-\kappa_1)\cdot ...\cdot
(t-\kappa_k)$ with $\kappa_j\in S^{UR}$, 
define the new polynomial $f\otimes g\in\C[t]$ with zeros
in $S^{UR}$ by
\begin{eqnarray}\label{2.4}
(f\otimes g)(t):=\prod_{i=1}^k\prod_{j=1}^l(t-\kappa_i\lambda_j).
\end{eqnarray}
Then
\begin{eqnarray}\label{2.5}
\divis (f\otimes g)&=& (\divis f)\cdot (\divis g),\\
\tr(\divis(f\otimes g))&=& (\tr(\divis f))\cdot (\tr(\divis g)),\\
\deg(f\otimes g)&=& (\deg f)\cdot (\deg g).\label{2.7}
\label{2.6}
\end{eqnarray}

The order $\ord(\zeta)\in\N$ of a unit root $\zeta\in S^{UR}$
is the minimal number $m\in\N$ with $\zeta^m=1$.
For $m\in \N$, the $m$-th cyclotomic polynomial is 
\begin{eqnarray}\label{2.8}
\Phi_m(t)=\prod_{\zeta:\, \ord(\zeta)=m}(t-\zeta)\in\C[t].
\end{eqnarray}
It is in $\Z[t]$, it has degree $\varphi(m)$, 
and it is irreducible in $\Z[t]$ and $\Q[t]$.
Denote 
\begin{eqnarray}\label{2.9}
\Lambda_m:=\divis(t^m-1),\quad \Psi_m:=\divis\Phi_m,\quad
E_m:=\frac{1}{m}\Lambda_m.
\end{eqnarray}
Then $\Lambda_1=\langle 1\rangle$. Of course
\begin{eqnarray}\label{2.10}
t^n-1&=&\prod_{m|n}\Phi_m(t),\hspace*{2cm} 
\Lambda_n=\sum_{m|n}\Psi_m,\\
\Phi_n&=&\prod_{m|n}(t^m-1)^{\mu_{Moeb}(\frac{n}{m})},\quad 
\Psi_n=\sum_{m|n}\mu_{Moeb}(\frac{n}{m})\cdot\Lambda_m.
\label{2.11}
\end{eqnarray}
Here $\mu_{Moeb}$ is the M\"obius function \cite{Ai}
\begin{eqnarray}\label{2.12}
\mu_{Moeb}:\N&\to&\{0,1,-1\},\\
m&\mapsto& \left\{\begin{array}{ll}
(-1)^r & \textup{if }m=p_1\cdot ...\cdot p_r\textup{ with }
p_1,...,p_r\\
 & \textup{different prime numbers},\\
0& \textup{else}  
\end{array}\right. \nonumber
\end{eqnarray}
(here $r=0$ is allowed, so $\mu_{Moeb}(1)=1$).
The traces of $\Lambda_m$ and $\Psi_m$ are 
\begin{eqnarray}\label{2.13}
\tr\Lambda_m&=&\left\{\begin{array}{ll}
1&\textup{if }m=1\\
0&\textup{if }m\geq 2,\end{array}\right. \\
\tr\Psi_m&=&\mu_{Moeb}(m).\label{2.14}
\end{eqnarray}
It is easy to see that 
\begin{eqnarray}\label{2.15}
\Lambda_a\cdot \Lambda_b &=&
\gcd(a,b)\cdot \Lambda_{\lcm(a,b)},\quad
E_a\cdot E_b=E_{\lcm(a,b)},\\
\langle \zeta\rangle\cdot \Lambda_b&=&\Lambda_b\qquad
\textup{if }\ord(\zeta)|b,\label{2.16}\\
\divis(f)\cdot\Lambda_b&=&
\deg f\cdot \Lambda_b\qquad\textup{if }f|(t^b-1).\label{2.17}
\end{eqnarray}
Especially
\begin{eqnarray}\label{2.18}
\Lambda_a\cdot \Lambda_b &=&
a\cdot \Lambda_b\qquad\textup{if }a|b.
\end{eqnarray}
It is more difficult to write down
formulas for $\Psi_a\cdot \Psi_b$. They can be cooked up
from the following special cases.
\begin{eqnarray}\label{2.19}
\Psi_a\cdot\Psi_b &=& \Psi_{a\cdot b}\qquad\textup{if }
\gcd(a,b)=1,\\
\Psi_{p^a}\cdot \Psi_{p^b} &=& \varphi(p^b)\cdot\Psi_{p^a}
\qquad\textup{if }p \textup{ is a prime number}\label{2.20}\\
&&\hspace*{2cm}\textup{and }a>b\geq 0,\nonumber \\
\Psi_{p^a}\cdot \Psi_{p^a} &=& \varphi(p^a)
\cdot\sum_{j=0}^a\Psi_{p^j}-p^{a-1}\cdot\Psi_{p^a}
\label{2.21}\\
&&\hspace*{2cm}\textup{if }p \textup{ is a prime number and }a>0.
\nonumber
\end{eqnarray}
Especially
\begin{eqnarray}\label{2.22}
\Psi_{2^a}\cdot \Psi_{2^a} &=& 2^{a-1}
\cdot\sum_{j=0}^{a-1}\Psi_{p^j}.
\end{eqnarray}

\medskip

Fix a finite set $M\subset \N$ and a map $\nu:\N\to\N_0$ with 
support $M$ (so $M=\{m\in\N\, |\, \nu(m)\neq 0\}$)
and define the unitary polynomial
\begin{eqnarray}\label{2.23}
\Delta:=\prod_{m\in M}\Phi_m^{\nu(m)}\in\Z[t]
\end{eqnarray}
(if $M=\emptyset$ then $\Delta=1$). Of course, then
\begin{eqnarray}\label{2.24}
\divis \Delta&=& \sum_{m\in M}\nu(m)\cdot\Psi_m.
\end{eqnarray}
Define also
\begin{eqnarray}\label{2.25}
d_M:=\lcm(m\in M\}.
\end{eqnarray}
Then $M\subset\{m\in\N\, |\, m|d_M\}$.
\eqref{2.11} and \eqref{2.24} give a unique function
$\chi:\N\to\Z$ with finite support 
\begin{eqnarray}\label{2.26}
\supp(\chi)\subset \{n\in\N\, |\, \exists\ m\in M
\textup{ with }n|m\}\subset\{n\in\N\, |\, n|d_M\}
\end{eqnarray}
and 
\begin{eqnarray}\label{2.27}
\divis \Delta &=& \sum_{n\in\N}\chi(n)\cdot \Lambda_n,\\
\nu(m)&=& \sum_{n:\, m|n}\chi(n),\label{2.28}\\
\chi(n)&=& \sum_{m:\, n|m} \nu(m)\cdot \mu_{Moeb}(\frac{m}{n}).
\label{2.29}
\end{eqnarray}
$\nu$ and $\chi$ and the following third function
$L:\N\to\Z$ determine each other. $L$ does not have finite support.
The numbers $L(k)\in\Z$ are the {\it Lefschetz numbers} of 
$\Delta$. They are defined by
\begin{eqnarray}\label{2.30}
L(k)&:=& \sum_{j=1}^{\deg \Delta}\lambda_j^k\qquad\textup{if }
\Delta=\prod_{j=1}^{\deg \Delta}(t-\lambda_j).
\end{eqnarray}
Especially $L(1)=\tr(\divis \Delta)$. Observe 
\begin{eqnarray}
\sum_{a=0}^{m-1}\langle e^{2\pi i a/m}\rangle&=& \Lambda_m,
\nonumber\\
\sum_{a=0}^{m-1}\langle e^{2\pi i ka/m}\rangle
&=&\gcd(k,m)\cdot \Lambda_{m/\gcd(k,m)}
\label{2.31}
\end{eqnarray}
Thus
\begin{eqnarray}
L(k)&=& \sum_{m\in M}\chi(m)\cdot\tr
\Bigl(\gcd(k,m)\Lambda_{m/\gcd(k,m)}\Bigr) \nonumber\\
&=& \sum_{m:\, m|k}m\cdot\chi(m).\label{2.32}
\end{eqnarray}
M\"obius inversion \cite{Ai} gives 
\begin{eqnarray}\label{2.33}
m\cdot\chi(m) &=& \sum_{k|m}\mu_{Moeb}(\frac{m}{k})\cdot L(k).
\end{eqnarray}
$L$ does not have finite support, but the following extended
periodicity property:
\begin{eqnarray}\label{2.34}
L(k)&=&L(\gcd(k,d_M)).
\end{eqnarray}
Therefore $L$ is determined by its values on 
$\{m\in\N\, |\, m|d_M\}$.
\eqref{2.34} implies the periodicity
\begin{eqnarray}\label{2.35}
L(k)=L(k+d_M),
\end{eqnarray}
but it is stronger. In fact,
\eqref{2.34} is equivalent to 
$\supp(\chi)\subset\{m\in\N\, |\, m|d_M\}$
and to $M\subset\{m\in \N\, |\, m|d_M\}$.

All the formulas \eqref{2.24} -- \eqref{2.35} make also sense
if $\divis\Delta$ is replaced by any element of $\Q\langle S^{UR}\rangle$.
Then $\nu,\chi$ and $L$ have values in $\Q$.

\section{Weight system and characteristic polynomial of a
quasihomogeneous singularity}\label{s3}
\setcounter{equation}{0}

\noindent
Fix a number $n\in \N$, and denote $N:=\{1,2,...,n\}$ and 
$e_i:=(\delta_{ij})_{j=1,...,n}\in\N_0^n$ for $i\in N$. 

\begin{definition}\label{t3.1}
A {\it weight system} is a tuple $(v_1,...,v_n,d)\subset(\Q_{>0})^{n+1}$
with $v_i<d$. Another weight system is {\it equivalent} to it,
if the second one has the form $q\cdot (v_1,...,v_n,d)$
for some $q\in\Q_{>0}$. 
A weight system is {\it integer} if $(v_1,...,v_n,d)\in\N^{n+1}$.
It is {\it reduced} if it is integer and it is minimal with this 
property, i.e. $\gcd(v_1,...,v_n,d)=1$. 
It is normalized if $d=1$. 
\end{definition}

Any equivalence class contains
a unique reduced weight system and a unique normalized weight system.
>From now on, the letters $(v_1,...,v_n,d)$ will be reserved
for integer weight systems, and $(w_1,...,w_n,1)$ will be the
equivalent normalized weight system, i.e. $w_i=\frac{v_i}{d}$.

Let $(v_1,...,v_n,d)$ be an integer weight system
(not necessarily reduced, it does not matter here).
For $J\subset N$ and $k\in\N_0$ denote
\begin{eqnarray*}
\Z^J&:=& \{\alpha\in\Z^n\ |\ \alpha_i=0\textup{ for }i\notin J\},
\qquad \N_0^J:=\Z^J\cap\N_0^n,\\
(\Z^n)_k&:=& \{\alpha\in \Z^n\ |\ \sum_i\alpha_i\cdot v_i=k\},
\qquad (\N_0^n)_k:=(\Z^n)_k\cap\N_0^n,\\
(\Z^J)_k&:=& \Z^J\cap (\Z^n)_k,
\qquad (\N_0^J)_k:=(\Z^J)_k\cap\N_0^n=\N_0^J\cap (\N_0^n)_k.
\end{eqnarray*}

\begin{remark}\label{t3.2}
For $J\subset N$ with $J\neq \emptyset$ define the semigroup
\begin{eqnarray}\label{3.1}
SG(J):=\sum_{j\in J}\N_0\cdot v_j\subset\N_0
\end{eqnarray}
and observe
\begin{eqnarray}\label{3.2}
\sum_{j\in J}\Z\cdot v_j=\Z\cdot \gcd(v_j\, |\, j\in J).
\end{eqnarray}
Then 
\begin{eqnarray}
(\N_0^J)_k\neq\emptyset&\iff& k\in SG(J),\label{3.3}\\
(\Z^J)_k\neq\emptyset&\iff& \gcd(v_j\, |\, j\in J) | k.
\label{3.4}
\end{eqnarray}
\end{remark}

The following combinatorial lemma is a specialization of  
\cite[Lemma 2.1]{HK1}.
It will be useful in theorem \ref{t3.5}.
(The conditions $(C2)'$ and $(C3)$ in \cite[Lemma 2.1]{HK1} 
are less important.)

\begin{lemma}\label{t3.3}
Fix an integer weight system $(v_1,...,v_n,d)$.
The following three conditions $(C1)$, $(C1)'$ and $(C2)$
are equivalent.
\begin{list}{}{}
\item[(C1):] 
\quad $\forall\ J\subset N\textup{ with }J\neq \emptyset
\quad (\N_0^J)_d\neq\emptyset$ or $\exists\ K\subset N-J$\\
\hspace*{3cm} with $|K|=|J|$ and 
$\forall\ k\in K\ (\N_0^J)_{d-v_k}\neq\emptyset.$\\
\item[(C1)':] 
\quad As (C1), but only for $J$ with $|J|\leq \frac{n+1}{2}$.\\
\item[(C2):] 
\quad $\forall\ J\subset N\textup{ with }J\neq \emptyset
\quad \exists\ K\subset N$ \\
\hspace*{3cm} with $|K|=|J|$ and $\forall\ k\in K\ 
(\N_0^J)_{d-v_k}\neq\emptyset$.
\end{list}
\end{lemma}

{\bf Proof:}
$(C1)\Rightarrow(C1)'$ is trivial.
$(C1)'\Rightarrow(C1)$ and $(C2)\Rightarrow(C1)$ are easy.
See \cite{HK1} for details. The 
least easy implication is $(C1)\Rightarrow(C2)$. 
In \cite{HK1} it was proved via the condition $(C3)$ there. 
A more direct proof will be given now.

Suppose that $(C1)$ holds.
Fix $J\subset N$ with $J\neq\emptyset$. 
We want to find a $K\subset N$ such that $J$ and $K$ 
satisfy $(C2)$.
Define the support of $J$ by 
\begin{eqnarray*}
\supp(J):=\{j\in J\, |\  \exists\ \alpha\in (\N_0^J)_d
\textup{ with }\alpha_j\neq 0\}\subset J.
\end{eqnarray*}
Consider $J_1:=J-\supp(J)$. 

{\bf 1st case,} $J_1=\emptyset$: Then $J$ and $K:=J$ satisfy
$(C2)$.

{\bf 2nd case,} $J_1\neq\emptyset$: The definition of $J_1$
implies $(\N_0^{J_1})_d=\emptyset$. Therefore $(C1)$ gives the
existence of a set $K_1\subset N-J_1$ with $|K_1|=|J_1|$ and
$\forall\ k\in K_1\quad (\N_0^{J_1})_{d-v_k}\neq\emptyset$.

Because of $v-d_k>0$, any element $\beta\in (\N_0^{J_1})_{d-v_k}$
satisfies $\beta_j\neq 0$ for some $j\in J_1$.

If some $k\in K_1$ would be in $J-J_1$, then
for any element $\beta\in (\N_0^{J_1})_{d-v_k}$ the element
$\alpha:=\beta+e_k$ would contradict 
$J_1\cap\supp(J)=\emptyset$.
Thus $K_1\subset N-J$. Now $J$ and $K:=K_1\cup\supp(J)$
satisfy $(C2)$.\hfill$\Box$

\begin{remarks}\label{t3.4}
(i) Denote with a bar the analogous conditions
$\oooo{(C1)}, \oooo{(C1)'}, \oooo{(C2)}$ 
where $\N_0$ is replaced by $\Z$. 
Also these conditions are equivalent to one another.
The proof is the same as above.

\medskip
(ii) Recall that a polynomial 
\begin{eqnarray*}
f=\sum_{\alpha\in\N_0^n}a_\alpha\cdot x^\alpha\in \C[x_1,...,x_n]
\quad\textup{where}\quad x^\alpha=x_1^{\alpha_1}...x_n^{\alpha_n}
\end{eqnarray*} 
is quasihomogeneous with respect to a weight system
$(v_1,...,v_n,d)$ if
\begin{eqnarray*}
\sum_{i=1}^n\alpha_i\cdot v_i=d\quad\textup{for all }\alpha
\textup{ with }a_\alpha\neq 0. 
\end{eqnarray*}
Recall that a quasihomogeneous polynomial has an isolated
singularity at 0 if the functions 
$\frac{\paa f}{\paa x_i}$ vanish simultaneously precisely at 0.
\end{remarks}

The following theorem is cited from \cite[Theorem 2.2]{HK1}.
It was first proved by Kouchnirenko \cite[Remarque 1.13 (i)]{Ko1}. 
See \cite[Remarks 2.3]{HK1}
for its history and contributions in \cite{Sa1} \cite{Ko1} 
\cite{Ko2} \cite{OR1} \cite{Sh} \cite{Wa} \cite{KS}.

\begin{theorem}\label{t3.5}
Let $(v_1,...,v_n,d)\in\N^{n+1}$ be an 
integer weight system.
The following conditions are equivalent.
\begin{list}{}{}
\item[(IS3):] \quad There exists a quasihomogeneous polynomial $f$ with
an isolated singularity at 0.
\item[(IS3)':] \quad A generic quasihomogeneous polynomial has an 
isolated singularity at 0.
\item[$(C1)$ to $(C2)$:] \quad $(v_1,...,v_n,d)$ satisfies 
$(C1),(C1)'$ and $(C2)$.
\end{list}
\end{theorem}

In definition \ref{t3.6}, some objects will be associated
to any weight system. Before studying them in the case of
weight systems of quasihomogeneous singularities, their
shape under weaker conditions will be discussed in lemma \ref{t3.7}.

\begin{definition}\label{t3.6}
Let $({\bf v},d)=(v_1,...,v_n,d)\in\N^{n+1}$ be an 
integer weight system. 
Define unique numbers $s_1,...,s_n,t_1,...,t_n\in\N$
by 
\begin{eqnarray*}
\frac{v_i}{d}=\frac{s_i}{t_i}\quad\textup{and}\quad
\gcd(s_i,t_i)=1.
\end{eqnarray*}
They depend only on the normalized weight system
${\bf w}=(w_1,...w_n)=(\frac{v_1}{d},...,\frac{v_n}{d})$.

Define 
\begin{eqnarray}\label{3.5}
d_{\bf w}:=\lcm(t_j\, |\,  j\in N).
\end{eqnarray}
Of course $d_{\bf w}|d$. If $({\bf v},d)$ is reduced
and $\gcd(v_1,...,v_n)|d$ (which holds for example
if $\oooo{(C2)}$ holds), then 
$\gcd(v_1,...,v_n)=1$ and then $d_{\bf w}=d$.

For $k\in\N$ define
\begin{eqnarray}\label{3.6}
M(k)&:=&\{j\in N\, |\, t_j|k\},\\
\textup{and }\mu(k)&:=&\prod_{j\in M(k)}(\frac{1}{w_j}-1)
=\prod_{j\in M(k)}\frac{d-v_j}{v_j}\in\Q_{>0}\label{3.7}
\end{eqnarray}
(the empty product is by definition $1$).
Define a quotient of polynomials
\begin{eqnarray}\label{3.8}
\rho_{({\bf v},d)}(t)&:=&t^{v_1+...+v_n}\cdot 
\prod_{j=1}^n \frac{t^{d-v_j}-1}{t^{v_j}-1}
\in\Q(t)
\end{eqnarray}
and an element of $\Q\langle S^{UR}\rangle$
\begin{eqnarray}\label{3.9}
D_{\bf w}&:=& \prod_{j=1}^n 
\left(\frac{1}{s_j}\Lambda_{t_j}-\Lambda_1\right) \in 
\Q\langle S^{UR} \rangle.
\end{eqnarray}
\end{definition}

\begin{lemma}\label{t3.7}
Let $({\bf v},d)$ be an integer weight system.

(a) Then 
\begin{eqnarray}\label{3.10}
M(k)&=&M(\gcd(k,d_{\bf w}))\\
&=& \{j\in N\, |\, \frac{d}{\gcd(k,d_{\bf w})} | v_j\},
\label{3.11}\\
\mu(k)&=&\mu(\gcd(k,d_{\bf w})).\label{3.12}
\end{eqnarray}

(b) The Lefschetz numbers $L(k)$ of the element $D_{\bf w}$ are
\begin{eqnarray}\label{3.13}
L(k)&=& L(\gcd(k,d_{\bf w}))\\
&=& (-1)^{n-|M(k)|}\cdot \mu(k)\in \Q^*.\label{3.14}
\end{eqnarray}

(c) 
\begin{eqnarray}\label{3.15}
({\bf v},d)\textup{ satisfies }\oooo{(C2)}&\Rightarrow& 
\mu(k)\in \N\quad\textup{for all }k\in\N.\\
({\bf v},d)\textup{ satisfies }\oooo{(C2)}&\iff& 
\rho_{({\bf v},d)}\in \Z[t].\label{3.16}
\end{eqnarray}
\end{lemma}

{\bf Proof:}
(a) As all $t_j$ divide $d_{\bf w}$ by definition of $d_{\bf w}$,
$t_j|k$ is equivalent to $t_j|\gcd(k,d_{\bf w})$.
This shows \eqref{3.10} and \eqref{3.12}.
Now suppose $k|d_{\bf w}$ (just for simplicity of notations).
By definition $t_j=d/\gcd(v_j,d)$. 
Thus for any $j\in J$
\begin{eqnarray*}
t_j|k
\iff \frac{d}{\gcd(v_j,d)}|k
\iff \frac{d}{k}| \gcd(v_j,d)
\iff \frac{d}{k}|v_j.
\end{eqnarray*}
This shows \eqref{3.11}.

(b) The following calculation gives \eqref{3.14}.
\begin{eqnarray*}
L(k)&=& \tr\left(\prod_{j=1}^n \left(\frac{\gcd(k,t_j)}{s_j}
\Lambda_{t_j/\gcd(k,t_j)}-\Lambda_1\right)\right)\\
&=& \prod_{j=1}^n \left( \frac{\gcd(k,t_j)}{s_j}\cdot
\tr\left(\Lambda_{t_j/\gcd(k,t_j)}\right)-1\right) \\
&=& \prod_{j\in M(k)}\left(\frac{t_j}{s_j}-1\right) 
\cdot \prod_{j\notin M(k)}(-1)\\
&=& (-1)^{n-|M(k)|}\cdot \mu(k).
\end{eqnarray*}
The equality \eqref{3.13} $L(k)=L(\gcd(k,d_{\bf w}))$ is a consequence
of the analogous properties \eqref{3.10}
of $M(k)$ and \eqref{3.12} of $\mu(k)$ and of \eqref{3.14}.

\medskip
(c) Recall from remark \ref{t3.2} 
that $(\Z^J)_{d-v_l}\neq \emptyset
\iff \gcd(v_j\, |\, j\in J)|(d-v_l)$. 
Therefore $\oooo{(C2)}$ says
\begin{list}{}{}
\item[$(GCD)$]
\quad $\forall \ J\subset N$ the $\gcd(v_j\, |\, j\in J)$ divides
at least $|J|$ of the numbers $d-v_l$ for $l\in N$.
\end{list}

For any $k\in \N$ with $k|d$, obviously the analogous
condition with $M(k)$ instead of $N$ holds then, too.
It is easy to derive from this directly $\mu(k)\in\N$ 
for this $k$.
But it will also follow from the consideration below of 
$\rho_{({\bf v},d)}$. 
Then $\mu(k)\in\N$ for all $k\in \N$ follows with \eqref{3.12}.


$\rho_{({\bf v},d)}$ is a quotient of cyclotomic polynomials.
The condition $(GCD)$ says that any cyclotomic polynomial in
the denominator turns up with at least the same multiplicity
in the numerator,
especially the cyclotomic polynomials
$\Phi_{\gcd(v_j\, |\, j\in J)}$ for some $J\subset N$. 
Thus $\rho_{({\bf v},d)}\in\Z[t]$ is equivalent
to $\oooo{(C2)}$.

Now suppose that $\oooo{(C2)}$ holds. 
Then $(GCD)$ holds for $N$ and also for any $M(k)$ instead of $N$.
The argument above for $\rho_{({\bf v},d)}\in\Z[t]$ 
applies also to the partial product
$\prod_{j\in M(k)}(...)$ for any $k\in \N$ and shows that it is 
in $\Z[t]$. Dividing out
all factors $(t-1)$, one can insert $t=1$ and obtains for the
partial product 
\begin{eqnarray*}
\mu_k=\prod_{j\in M(k)}\frac{d-v_j}{v_j}\in\Z\cap \Q_{>0}=\N.
\end{eqnarray*}
\hfill$\Box$ 

\begin{remarks}\label{t3.8}
(i) Let $({\bf v},d)$ be an integer weight system with
$\oooo{(C2)}$. Then \eqref{3.16} gives 
$\rho_{({\bf v},d)}=\sum_{\alpha\in\frac{1}{d}\Z}
\sigma(\alpha)\cdot t^{d\cdot\alpha}$ with a function
$\sigma:\frac{1}{d}\Z\to\Z$ with finite support.
And \eqref{3.15} gives $L(k)\in\Z$.

{\bf Open problem 1:}
\begin{list}{}{}
\item[(a)]
Is $D_{\bf w}\in\Z\langle S^{UR}\rangle?$
Equivalent: Are the $\chi(m)$ which are determined by the
$L(k)$ and \eqref{2.33} in $\Z$?
\item[(b)]
Is $D_{\bf w}=\sum_{\alpha\in\frac{1}{d}\Z}
\sigma(\alpha)\cdot \langle e^{2\pi i\alpha}\rangle$?
\end{list}

{\it Yes} for problem 1 (b) would imply {\it Yes} for 
problem 1 (a). 
\end{remarks}

The following theorem is classical, see the remarks \ref{t3.10}
for its origins.

\begin{theorem}\label{t3.9}
Let $({\bf v},d)$ be an integer weight system with $(C2)$,
i.e. a weight system of isolated quasihomogeneous singularities.

Then the divisor of the characteristic polynomial 
of its monodromy is $D_{\bf w}$,
so here $D_{\bf w}\in \N_0[t]$, and 
\begin{eqnarray}\label{3.19}
\tr((\textup{monodromy})^k)=(-1)^{n-|M(k)|}\cdot \mu(k).
\end{eqnarray}
Also $\rho_{({\bf v},d)}\in\N_0[t]$, thus
\begin{eqnarray}\label{3.20}
\rho_{({\bf v},d)}=\sum_{j=1}^\mu t^{d\cdot\alpha_j}\quad 
\textup{for certain }\alpha_j\in\Q.
\end{eqnarray}
These numbers $(\alpha_1,...,\alpha_\mu)$ are the 
{\rm exponents} of the singularity,
and $e^{2\pi i\alpha_1},...,e^{2\pi i\alpha_\mu}$
are the eigenvalues of the monodromy, i.e. the zeros of 
the characteristic polynomial,
so here
\begin{eqnarray}\label{3.21}
D_{\bf w}&=& \sum_{j=1}^\mu \langle e^{2\pi i\alpha_j}\rangle
\in \N_0\langle S^{UR}\rangle.
\end{eqnarray}
\end{theorem}

\begin{remarks}\label{t3.10}
(i) Formula \eqref{3.19} was shown by Milnor in \cite[\S 9.6]{Mi}. 
Of course, the trace of the $k$-th power of the monodromy is
precisely the $k$-th Lefschetz number of the characteristic
polynomial of the monodromy. Therefore \eqref{3.19}
together with \eqref{3.14} and 
the equivalence of the data $L,\chi,\nu,\Delta$ 
in section \ref{s2}
implies that the characteristic polynomial of the monodromy has
the divisor $D_{\bf w}$. This was first seen in \cite{MO}.

\medskip
(ii)
The polynomial $\rho_{({\bf v},d)}$ is the generating function
of the exponents, which are up to the shift $v_1+...+v_n$ 
the weighted degrees of the Jacobi algebra
\begin{eqnarray*}
\frac{\C\{x_1,...,x_n\}}{\left(\frac{\paa f}{\paa x_1}, ... ,
\frac{\paa f}{\paa x_n}\right)} \stackrel{f\textup{ q.h.}}{\cong}
\frac{\C[x_1,...,x_n]}{\left(\frac{\paa f}{\paa x_1},...,
\frac{\paa f}{\paa x_n}\right)}.
\end{eqnarray*}
This was (re)discovered by many people.
Therefore $\rho_{({\bf v},d)}\in\N_0[t]$.

\medskip
(iii) Let $({\bf v},d)$ be an integer weight system with
$n\leq 3$. Theorem 3 in \cite{Sa2} says
\begin{eqnarray*}
\rho_{({\bf v},d)}\in\Z[t]\iff (IS3).
\end{eqnarray*}
With \eqref{3.16} and theorem \ref{t3.5}, this is equivalent
to $\oooo{(C1)}\iff (C1)$ for $n\leq 3$.
This equivalence $\oooo{(C1)}\iff (C1)$ for $n\leq 3$
is lemma 2.5 in \cite{HK1}. It has a short combinatorial proof.
\end{remarks}

\begin{remarks}\label{t3.11}
(i) Theorem \ref{t3.9} implies that for an integer weight
system $({\bf v},d)$ with $(C2)$ the answer to the
problem 1 (a)+(b) is {\it Yes}.
But the proof is not combinatorial.
Theorem \ref{t3.9} gives also the following three 
implications:
\begin{eqnarray}\label{3.22}
(C2)&\Rightarrow& D_{\bf w}\in \N_0\langle S^{UR}\rangle,\\
\quad (C2)&\Rightarrow& \rho_{({\bf v},d)}\in \N_0[t],
\label{3.23}\\
(C2)&\Rightarrow& D_{\bf w}=\sum_{j=1}^\mu 
\langle e^{2\pi i\alpha_j}\rangle .\label{3.24}
\end{eqnarray}
\eqref{3.24} is the positive answer to problem 1 (b)
in the case $(C2)$. The known proofs of \eqref{3.22}, 
\eqref{3.23} and \eqref{3.24} are not combinatorial.

\medskip
(ii) Let $({\bf v},d)$ be an integer weight system with
$(C2)$. 

{\bf Open problem 2:} 
\begin{list}{}{}
\item[(a)] Give a combinatorial proof of \eqref{3.22}.
\item[(b)] Give a combinatorial proof of 
\eqref{3.23}.
\item[(c)] Give a combinatorial proof of 
\eqref{3.24}.
\end{list}

\medskip
(iii) One can separate in $(C2)$ and $(C1)$ the conditions for $J$
of different values of $|J|\in N$. The conditions for $J$ with
$|J|=1$ lead to the graphs and types of a quasihomogeneous
singularity which are discussed in \cite[ch. 3]{HK1}.
The sections 4 and 6 in \cite{HK1} make extensive and
successful use of the conditions for $|J|=1$. 
Below in section \ref{s4}, we will extend formulas in \cite{HK1} 
for parts of the Milnor number $\mu$ 
to formulas for parts of $D_{\bf w}$.

But it is irritatingly difficult to make use of the conditions
for $J$ with $|J|\geq 2$ in $(C2)$ or $(C1)$. 
Though they must be used in solutions
of the problems in (iii), and probably also in a positive solution
of the conjecture \ref{t6.3} 
in section \ref{s6}, if that has a positive solution.

{\bf Open problem 3:} Make some good use of the conditions for $J$ 
with $|J|\geq 2$ in $(C2)$ or $(C1)$.
\end{remarks}

The last point in this section is a discussion of a well known
fact on the order of the monodromy of a quasihomogeneous
singularity. That order is 
\begin{eqnarray*}
d_{mon}:=\lcm(m\in\N\, |\, \nu(m)>0)\qquad \textup{where }
D_{\bf w}=\sum_{m\in \N}\nu(m)\cdot \Psi_m.
\end{eqnarray*}

\begin{lemma}\label{t3.12}
In the case of a quasihomogeneous singularity 
with $w_j\leq\frac{1}{2}$ for all $j\in N$, 
$d_{mon}=d_{\bf w}$ or $d_{mon}=\frac{d_{\bf w}}{2}$.
If all $w_j<\frac{1}{2}$ then $d_{mon}=d_{\bf w}$.
\end{lemma}

{\bf Proof:}
Because of the definition \eqref{3.9} of $D_{\bf w}$, 
$d_{mon}$ is a divisor of $d_{\bf w}$.
The equalities
\begin{eqnarray*}
\prod_{j=1}^n\left(\frac{1}{w_j}-1\right) 
&=&\mu = \tr \id = \tr (Mon)^{d_{mon}}
=L(d_{mon})\\
&=& \pm\mu(d_{mon})=\pm
\prod_{j\, :\, t_j|d_{mon}}\left(\frac{1}{w_j}-1\right)
\end{eqnarray*}
show that the second product can miss only indices $j$ 
with $w_j=\frac{1}{2}$.
Therefore $\lcm(t_j\ |\ w_j<\frac{1}{2})$ divides $d_{mon}$. 
If all $w_j<\frac{1}{2}$ then 
$\lcm(t_j\ |\ w_j<\frac{1}{2})=d_{\bf w}$.
If some $w_j=\frac{1}{2}$ then 
$\lcm(t_j\ |\ w_j<\frac{1}{2})=d_{\bf w}$ or $\frac{d_{\bf w}}{2}$.
\hfill$\Box$

\section{Formulas for quasihomogeneous singularities
of cycle type and of chain type}\label{s4}
\setcounter{equation}{0}

\noindent
The formulas in this section concern the quasihomogeneous
singularities of cycle type and of chain type. 

They start with normalized weight systems 
$({\bf w},1)=(w_1,...,w_n,1)$ which satisfy a part of the
conditions $(C1)$ and $(C2)$ in lemma \ref{t3.3}, 
namely that part which concerns subsets $J\subset N$ with $|J|=1$.
That part leads to graphs and types of weight systems,
see section 3 in \cite{HK1}. As already said in remark \ref{t3.10}
(iv), it is difficult to make use of the conditions in
$(C1)$ and $(C2)$ for $J$ with $|J|\geq 2$. The formulas here 
do not make use of these higher conditions.

The formulas extend formulas in \cite{HK1} for parts of
the weight system and parts of the Milnor number to formulas for
parts of $D_{\bf w}$. Some calculations already made
in \cite{HK1} will be reproduced here for better readability.
We start with the cycle type,
then consider a generalization of the chain type and
finally specialize that to the chain type.
The formulas for the generalization of the chain type will
allow to glue its root to another graph.

Define the function
\begin{eqnarray}\label{4.1}
\rho:\bigcup_{k=0}^\infty\Z^k&\to&\Z,\\
\rho(x_1,...,x_k)&:=& x_1...x_k-x_2...x_k+...+(-1)^{k-1}x_k+(-1)^k
\nonumber
\end{eqnarray}
(the case $k=0$ is $\rho(\emptyset)=1$).

\begin{lemma}\label{t4.1} (Partly \cite[Lemma 3.4 and (4.6)]{HK1})
Fix $n\in\N$ and $n$ numbers $a_1,...,a_n\in\N$ with
$\prod_{j=1}^n a_j-(-1)^n>0$ and, if $n$ is even,
neither $a_j=1$ for all even $j$ nor $a_j=1$ for all odd $j$.

Then there is a unique normalized weight system
$({\bf w},1)=(w_1,...,w_n,1)$ with 
$a_jw_j+w_{j-1}=1$ for all $j\in N$, where $w_0:=w_n$. It is
\begin{eqnarray}\label{4.2}
w_j&=& \frac{v_j}{d}\qquad\textup{where }\\
v_j&:=&\rho(a_{j-1},a_{j-2},...,a_2,a_1,a_n,a_{n-1},...,a_{j+1}),
\nonumber\\
d&:=&\prod_{j=1}^na_j-(-1)^n.\nonumber
\end{eqnarray}
Define
\begin{eqnarray*}
\gamma:=\gcd(v_1,d).
\end{eqnarray*}
Then the unique numbers $s_j,t_j\in\N$ with $\gcd(s_j,t_j)=1$
and $w_j=\frac{s_j}{t_j}$ from definition \ref{t3.6} are
\begin{eqnarray}\label{4.3}
s_j&=& \frac{v_j}{\gamma},\\
t_1&=&...=t_n=\frac{d}{\gamma}.\nonumber
\end{eqnarray}
Especially $\gamma=\gcd(v_j,d)$ for any $j\in N$.

The Milnor number is
\begin{eqnarray}\label{4.4}
\mu&=& \prod_{j=1}^n a_j.
\end{eqnarray}
The divisor $D_{\bf w}$ from definition \ref{t3.6} is
\begin{eqnarray}\label{4.5}
D_{\bf w}&=& 
\gamma\cdot\Lambda_{d/\gamma} + (-1)^n\cdot\Lambda_1.
\end{eqnarray}
\end{lemma}

{\bf Proof:} 
The matrix of the system $a_jw_j+w_{j-1}=1$ of linear equations has
the determinant 
\begin{eqnarray*}
\det\begin{pmatrix}a_1 & & & 1\\ 1 & a_2 & & \\ 
 & \ddots & \ddots &  \\  & & 1 & a_n \end{pmatrix}
=\prod_{j=1}^n a_j - (-1)^n=d>0,
\end{eqnarray*}
here $d>0$ by hypothesis. Therefore it has a unique solution.
It is easy to see that this solution is given by \eqref{4.2}.
The conditions that in the case $n$ even neither $a_j=1$ for
all even $j$ nor $a_j=1$ for all odd $j$ make sure that
the numbers $v_j$ and the weights $w_j$ are not zero, but positive.
The equation $a_jw_j+w_{j-1}=1$ implies $w_j<1$.

By definition $t_1=d/\gamma$. The identities
(where $w_0=w_n,s_0=s_n,t_0=t_n$)
\begin{eqnarray}\label{4.6}
\left. \begin{array}{r}
\frac{s_j}{t_j}=w_j=\frac{1-w_{j-1}}{a_j}
=\frac{t_{j-1}-s_{j-1}}{t_{j-1}\cdot a_j}\\
\gcd(t_{j-1},t_{j-1}-s_{j-1})=1
\end{array} \right\}\Rightarrow 
t_j=t_{j-1}\cdot\frac{a_j}{\gcd(a_j,t_{j-1}-s_{j-1})}
\end{eqnarray}
show $t_{j-1}|t_j$. As we have a cycle here,
$t_j=d/\gamma$ and $\gcd(v_j,d)=\gamma$ for any $j\in N$.

The Milnor number is calculated by (with $v_0=v_n$)
\begin{eqnarray*}\nonumber
\mu= \prod_{j=1}^n\frac{d-v_{j-1}}{v_{j-1}} 
=\prod_{j=1}^n \frac{a_j\cdot v_j}{v_{j-1}}
=a_1\cdot ...\cdot a_n.
\end{eqnarray*}
The divisor $D_{\bf w}$ is defined in \eqref{3.9}.
Because of \eqref{2.15} it has only the two summands
$\Lambda_{d/\gamma}$ and $\Lambda_1$, and the coefficient
of $\Lambda_1$ is obviously (from \eqref{3.9}) $\chi(1)=(-1)^n$.
As $\mu=\deg D_{\bf w}= 
\chi(d/\gamma)\cdot d/\gamma + \chi(1)$, the 
coefficient $\chi(d/\gamma)$ of $\Lambda_{d/\gamma}$ is
$\chi(d/\gamma)=\gamma$, so \eqref{4.5} holds.
\hfill$\Box$

\begin{lemma}\label{t4.2}
(Partly \cite[(4.10)]{HK1})
Fix $n\in\N$, $n$ numbers $a_1,...,a_n\in\N$, two numbers
$s_0,t_0\in\N$ with $s_0<t_0$ and $\gcd(s_0,t_0)=1$,
and define $w_0:=\frac{s_0}{t_0}$. 

Then there are unique weights $w_1,...,w_n\in\Q\cap (0,1)$
with $a_jw_j+w_{j-1}=1$ for $j=1,...,n$. 
Write $w_j=\frac{s_j}{t_j}$ with $s_j,t_j\in \N$ and
$\gcd(s_j,t_j)=1$ and 
$\beta_j:=\gcd(t_{j-1}-s_{j-1},a_j)\in\N$ and 
$\alpha_j:=\frac{a_j}{\beta_j}\in\N$. Then

\begin{eqnarray}\label{4.7}
s_j&=& \frac{\rho(a_{j-1},...,a_1)\cdot t_0+(-1)^js_0}
{\beta_j\cdot
...\cdot \beta_1},\\
t_j&=& \alpha_j\cdot t_{j-1}= 
\alpha_j\cdot ...\cdot \alpha_1\cdot t_0.\label{4.8}
\end{eqnarray}
The partial divisor and the partial Milnor number 
associated to $(w_1,...,w_n)$ are
\begin{eqnarray}\label{4.9}
\prod_{j=1}^n\left(\frac{1}{s_j}\Lambda_{t_j}-\Lambda_1\right)
&=& (-1)^n\Lambda_1+ \sum_{j=1}^n
(-1)^{n-j}\frac{\beta_j...\beta_1}{t_0-s_0}\cdot\Lambda_{t_j}\\
&=& (-1)^n E_1+\sum_{j=1}^n 
(-1)^{n-j}\frac{a_j...a_1}{1-w_0}
\cdot E_{t_j},\label{4.10}\\
\prod_{j=1}^n\left(\frac{1}{w_j}-1\right)
&=& \frac{\rho(a_n,a_{n-1},...,a_1)+(-1)^{n-1}w_0}{1-w_0}.
\label{4.11}
\end{eqnarray}
\end{lemma}

{\bf Proof:}
The weights are unique and in $\Q\cap(0,1)$ because
they are determined inductively by the equations 
$a_jw_j+w_{j-1}=1$, i.e.
\begin{eqnarray*}
w_j=\frac{1-w_{j-1}}{a_j}
= \frac{t_{j-1}-s_{j-1}}{a_j\cdot t_{j-1}}\\
= \frac{(t_{j-1}-s_{j-1})/\beta_j}{\alpha_j\cdot t_{j-1}}.
\end{eqnarray*}
As $1=\gcd(s_{j-1},t_{j-1})=\gcd(t_{j-1}-s_{j-1},t_{j-1})$,
this shows \eqref{4.8}. 
For $j=1$ \eqref{4.7} is clear. For $j\geq 2$ 
the additional calculation
\begin{eqnarray*}
\beta_j\cdot s_j&=& t_{j-1}-s_{j-1}\\
&=& \frac{a_{j-1}\cdot ...\cdot a_1\cdot t_0}{\beta_{j-1}\cdot
...\cdot \beta_1}-
\frac{\rho(a_{j-2},...,a_1)\cdot t_0+(-1)^{j-1}s_0}
{\beta_{j-1}\cdot ...\cdot \beta_1}\\
&=& \frac{\rho(a_{j-1},...,a_1)\cdot t_0+(-1)^js_0}
{\beta_{j-1}\cdot ...\cdot \beta_1}
\end{eqnarray*}
shows \eqref{4.7}.
Now the partial divisor is also calculated inductively.
The induction uses the partial divisor and the
partial Milnor number for $n-1$.
Also $t_j|t_n$ and \eqref{2.18} 
($\Lambda_a\Lambda_b=a\Lambda_b$ for $a|b$) are used. 
\begin{eqnarray*}
&&\left(\prod_{j=1}^{n-1}\left(\frac{1}{s_j}\Lambda_{t_j}-\Lambda_1\right)\right)
\cdot \left(\frac{1}{s_n}\Lambda_{t_n}-\Lambda_1\right)\\
&=& \left(\prod_{j=1}^{n-1}
\left(\frac{1}{w_j}-1\right)\right)\cdot 
\frac{1}{s_n}\Lambda_{t_n}
-\prod_{j=1}^{n-1}
\left(\frac{1}{s_j}\Lambda_{t_j}-\Lambda_1\right)\\
&=& \frac{\beta_n...\beta_1}{t_0-s_0}\cdot\Lambda_{t_n}+
(-1)^n+ \sum_{j=1}^{n-1}
(-1)^{n-j}\frac{\beta_j...\beta_1}{t_0-s_0}\cdot\Lambda_{t_j}.
\end{eqnarray*}
This shows \eqref{4.9} and \eqref{4.10}.
The partial Milnor number is the degree of the partial divisor.
\hfill$\Box$ 

\bigskip

Lemma \ref{t4.2} is a slight generalization of the
chain type. The following corollary specializes it
to the chain type.

\begin{corollary}\label{t4.3}
In the situation of lemma \ref{t4.2}, suppose $w_0=w_1$.
Then 
\begin{eqnarray}\label{4.12}
s_0=s_1=1,t_0=t_1=a_1+1,\beta_1=a_1,\alpha_1=1.
\end{eqnarray}
Define 
\begin{eqnarray}\label{4.13}
b_k&:=& (a_1+1)\cdot a_2\cdot ...\cdot a_k \quad
\textup{for }k=1,...,n,\quad b_0:= 1,\\
\mu_k&:=& \rho(a_k,...,a_2,a_1+1),\quad
\textup{for }k=1,...,n,\quad \mu_0:=1.\label{4.14}
\end{eqnarray}
Then
\begin{eqnarray}\label{4.15}
s_j&=& \frac{\mu_{j-1}}
{\beta_j\cdot
...\cdot \beta_2},\\
t_j&=& \alpha_j\cdot t_{j-1}= 
\alpha_j\cdot ...\cdot \alpha_2\cdot (a_1+1)
=\frac{b_j}
{\beta_j\cdot ...\cdot\beta_2},\hspace*{1cm}\label{4.16}
\end{eqnarray}
\begin{eqnarray}
D_{\bf w}&=&
\prod_{j=1}^n\left(\frac{1}{s_j}\Lambda_{t_j}-\Lambda_1\right)
\nonumber\\
&=& (-1)^n+ \sum_{j=1}^n
(-1)^{n-j}\beta_j...\beta_2\cdot\Lambda_{t_j}\label{4.17}\\
&=& (-1)^n+\sum_{j=1}^n 
(-1)^{n-j}b_j
\cdot E_{t_j},\hspace*{1cm}\label{4.18}\\
\mu_k&=& \prod_{j=1}^k\left(\frac{1}{w_j}-1\right)
=b_k-\mu_{k-1} .
\label{4.19}
\end{eqnarray}
Furthermore, define
\begin{eqnarray}\label{4.20}
\sum_{j=0}^n(-1)^{n-j}\Lambda_{b_j}&=:&
\sum_{i=1}^\mu\langle\lambda_i\rangle .
\end{eqnarray}
The definition \eqref{4.20} makes sense, as obviously
for the divisor on the left hand side,
\begin{eqnarray}\label{4.21}
\nu(m)=\left\{\begin{array}{ll}
1& \textup{if for some even }k\ m|b_{n-k},m\not|b_{n-k-1},\\
0& \textup{else.}\end{array}\right. 
\end{eqnarray}
Then
\begin{eqnarray}\label{4.22}
\sum_{j=0}^n(-1)^{n-j}\Lambda_{b_j}
&=& 
\prod_{j=1}^n\left(\frac{1}{\mu_{j-1}}
\Lambda_{b_j}-\Lambda_1\right)\\
D_{\bf w}
&=& \sum_{i=1}^\mu \langle \lambda_i^\mu\rangle .
\label{4.23}
\end{eqnarray}
\end{corollary}

{\bf Proof:}
Formula \eqref{4.12} is trivial.
The formulas \eqref{4.15} to \eqref{4.19} are immediate
consequences of the formulas in lemma \ref{t4.2}.
\eqref{4.22} is proved inductively by a similar
calculation as \eqref{4.9},
\begin{eqnarray*}
&&\left(\prod_{j=1}^n\left(\frac{1}{\mu_{j-1}}\Lambda_{b_j}-\Lambda_1\right)\right)
\cdot \left(\frac{1}{\mu_{n-1}}\Lambda_{b_n}-\Lambda_1\right)\\
&=& \left(\prod_{j=1}^{n-1}
\left(\frac{1}{w_j}-1\right)\right)\cdot 
\frac{1}{\mu_ {n-1}}\Lambda_{b_n}
-\prod_{j=1}^{n-1}
\left(\frac{1}{\mu_{j-1}}\Lambda_{b_j}-\Lambda_1\right)\\
&=& \Lambda_{b_n}+ \sum_{j=0}^{n-1}
(-1)^{n-j}\Lambda_{b_j}.
\end{eqnarray*}
For the final formula \eqref{4.23}, it is in view of
\eqref{2.31} and \eqref{4.18} enough to show
\begin{eqnarray}\label{4.24}
\frac{b_j}{\gcd(b_j,\mu)}&=& t_j\qquad \textup{for }j\geq 1.
\end{eqnarray}
But $\mu_k=b_k-\mu_{k-1}$ (for $k\geq 1$) and 
$b_k=a_kb_{k-1}$ (for $k\geq 2$) show
\begin{eqnarray*}
\gcd(b_j,\mu)&=& \gcd(b_j,\mu_{n-1})=...=
\gcd(b_j,\mu_j)=\gcd(b_j,\mu_{j-1}).
\end{eqnarray*}
As $w_j=\frac{s_j}{t_j}=\frac{\mu_{j-1}}{b_j}$, 
\begin{eqnarray*}
t_j = \frac{b_j}{\gcd(b_j,\mu_{j-1})}
=\frac{b_j}{\gcd(b_j,\mu)}.
\end{eqnarray*}
\hfill$\Box$

\begin{remark}\label{t4.4}
The last formula \eqref{4.23} in corollary \ref{t4.3}
fits to a result of Orlik and Randell \cite[(2.11) theorem]{OR2}.
They showed that the integral monodromy is the $\mu$-th 
power of a cyclic automorphism of the Milnor lattice,
whose eigenvalues are given by the divisor in \eqref{4.20}.
Formula \eqref{4.23} just confirms that the divisor
$D_{\bf w}$ has the eigenvalues which fit to this theorem.

We made this calculation mainly to see how it works
and to get some inspiration for good guesses for 
other types of weight systems of quasihomogeneous singularities.
\end{remark}

\section{Examples and counterexamples}\label{s5}
\setcounter{equation}{0}

\noindent
This section offers examples. Some of them are 
counterexamples to conjectures or hopes. 

\begin{example}\label{t5.1}
We begin with an example of Ivlev \cite[12.3]{AGV}.
It is the integer weight system 
\begin{eqnarray}\label{5.1}
({\bf v},d)=(1,24,33,58,265).
\end{eqnarray}
It satisfies $\oooo{(C1)}$, but not $(C1)$.
Observe
\begin{eqnarray*}
({\bf w},1)=(\frac{1}{265},\frac{24}{265},\frac{33}{265},
\frac{58}{265},1),\quad
w_j=\frac{v_j}{d}=\frac{s_j}{t_j}, 
\textup{ with }s_j=v_j,t_j=d,
\end{eqnarray*}
and 
\begin{eqnarray*}
265=5\cdot 53,\quad 264=3\cdot 8\cdot 11 = 8\cdot 33
=11\cdot 24,\\
265-33=232=4\cdot 58,\\
\gcd(24,33)=3, 265-58=207=3\cdot 69,\\
\textup{but }207\notin SG(24,33):=\N_0\cdot 24+\N_0\cdot 33.
\end{eqnarray*}
The following table lists the sets $J$ with $|J|\leq 2$ 
which satisfy
alone or with a suitable set $K\subset N-J$ the condition $(C1)$.
\begin{eqnarray*}
\begin{array}{l|l|l|l|l|l|l|l|l|l}
J&\{1\} & \{2\} & \{3\} & \{4\} & \{1,2\} & \{1,3\} & \{1,4\} & \{2,4\} & \{3,4\} \\ \hline
K&  & \{1\} & \{1\} & \{3\} &  &  &  & \{1,3\} &  
\end{array}
\end{eqnarray*}
The set $J=\{2,3\}$ satisfies with $K=\{1,4\}$ $\oooo{(C1)}$, 
but not $(C1)$.
In the notation of \cite[Example 3.2 (iii)]{HK1}, 
the weight system is of
type XII (but with a different numbering).
Ivlev (and we, too) calculated that 
$\rho_{({\bf v},d)}\in\N_0[t]$
(and not only in $\Z[t]$). 

The sets $M(k)$ are 
\begin{eqnarray*}
M(k) &=& 
\left\{ \begin{array}{ll}
M(265)=N=\{1,2,3,4\} & \textup{if }265|k,\\
M(1)=\emptyset & \textup{if }265\not| k. \end{array}\right. 
\end{eqnarray*}
Therefore only the values of $L(k)=(-1)^{n-|M((k)|}\cdot \mu(k)$
and $\chi(k)$ for $k\in\{1,265\}$ are interesting. 
\begin{eqnarray*}
(L(265),L(1)) &=& (66516,1)=(\mu,1),\\
(\chi(265),\chi(1)) &=& (251,1),\\
D_{\bf w} &=& 251\cdot \Lambda_{265}+1\cdot\Lambda_1 
=\frac{\mu-1}{265}\cdot \Lambda_{265}+\Lambda_1.
\end{eqnarray*}
\end{example}

\begin{remarks}\label{t5.2}
{\bf Open problem 4:} 
Find other examples of integer weight systems $({\bf v},d)$ 
which satisfy $\oooo{(C1)}$, but not $(C1)$.
Both cases, $\rho_{({\bf v},d)}\in\N_0[t]$
and $\rho_{({\bf v},d)}\in\Z[t]-\N_0[t]$, are interesting.
Because of remark \ref{t3.10} (iii), 
all such example satisfy $n\geq 4$.
Find examples with $n=4$ of other types as Ivlev's example,
which is of type XII in the notation of 
\cite[Example 3.2 (iii)]{HK1}.
\end{remarks}

\begin{examples}\label{t5.3}
Here some examples of weight systems of 
quasihomogeneous singularities are
given, together with the values of $\nu,\chi$ and $L$
from section \ref{s2}.

\medskip
(i) $n=3,\ N=\{1,2,3\}$,
\begin{eqnarray*}
(w_1,w_2,w_3,1)=(\frac{1}{4},\frac{1}{6},\frac{5}{12},1).
\end{eqnarray*}
The monomials $x_1^4,x_2^6,x_2x_3^2$ give the type II
in \cite[example 3.2 (ii)]{HK1}.
The following table lists all sets $M(k)$ and suitable 
values of $k$.
\begin{eqnarray*}
\begin{array}{cccc}
N&\{1\}&\{2\}&\emptyset\\ \hline
12 & 4 & 6 & 1
\end{array}
\end{eqnarray*}
Therefore only the values of $L(k)=(-1)^{n-|M(k)|}\cdot \mu(k)$
and $\chi(k)$ for $k\in\{12,4,6,1\}$ are interesting.
\begin{eqnarray*}
(L(12),L(4),L(6),L(1))=(21,3,5,-1),\\
(\chi(12),\chi(4),\chi(6),\chi(1))=(1,1,1,-1),\\
D_{\bf w}=\Lambda_{12}+\Lambda_4+\Lambda_6-\Lambda_1\\
= \Lambda_{12}+(\Psi_6+\Psi_4+\Psi_3+\Psi_2+\Psi_1)+\Psi_2.
\end{eqnarray*}

\medskip
(ii) $n=4,\ N=\{1,2,3,4\}$,
\begin{eqnarray*}
(w_1,w_2,w_3,w_4,1)=(\frac{1}{5},\frac{2}{5},\frac{1}{6},
\frac{5}{12},1).
\end{eqnarray*}
The monomials $x_1^5,x_1x_2^2,x_3^6,x_3x_4^2$ give the type 
XIII in \cite[example 3.2 (iii)]{HK1}.
The following table lists all sets $M(k)$ and suitable 
values of $k$.
\begin{eqnarray*}
\begin{array}{cccccc}
N&\{1,2,3\}&\{3,4\}&\{1,2\}&\{3\}&\emptyset\\ \hline
60 & 30 & 12 & 5 & 6 & 1
\end{array}
\end{eqnarray*}
Therefore only the values of $L(k)=(-1)^{n-|M(k)|}\cdot \mu(k)$
and $\chi(k)$ for $k\in\{60,30,12,5,6,1\}$ are interesting.
\begin{eqnarray*}
(L(60),L(30),L(12),L(5),L(6),L(1))=(42,-30,7,6,-5,1),\\
(\chi(60),\chi(30),\chi(12),\chi(5),\chi(6),\chi(1))=
(1,-1,1,1,-1,1),\\
D_{\bf w}=\Lambda_{60}-\Lambda_{30}+\Lambda_{12}+\Lambda_5
- \Lambda_6+\Lambda_1\\
= (\Psi_{60}+\Psi_{20}+\Psi_{12}+\Psi_5+\Psi_4+\Psi_1)+
(\Psi_{12}+\Psi_4+\Psi_1).
\end{eqnarray*}

\medskip
(iii) The curve singularity $D_{2q}$: 
\begin{eqnarray*}
n=2,N=\{1,2\},\mu=2q, (w_1,w_2,1)=(\frac{1}{2q-1},
\frac{q-1}{2q-1},1).
\end{eqnarray*}
The monomials $x_1^{2q-1},x_1x_2^2$ give the type II
in \cite[example 3.2 (i)]{HK1}.
The following table lists all sets $M(k)$ and suitable 
values of $k$.
\begin{eqnarray*}
\begin{array}{cc}
N&\emptyset\\ \hline
2q-1  & 1
\end{array}
\end{eqnarray*}
Therefore only the values of $L(k)=(-1)^{n-|M(k)|}\cdot \mu(k)$
and $\chi(k)$ for $k\in\{2q-1,1\}$ are interesting.
\begin{eqnarray*}
(L(2q-1),L(1))=(2q,1),\\
(\chi(2q-1),\chi(1))=(1,1),\\
D_{\bf w}=\Lambda_{2q-1}+\Lambda_1.
\end{eqnarray*}

\medskip
(iv) The curve singularity $D_{2q+1}$: 
\begin{eqnarray*}
n=2,N=\{1,2\},\mu=2q+1, (w_1,w_2,1)=(\frac{1}{2q},
\frac{2q-1}{4q},1).
\end{eqnarray*}
The monomials $x_1^{2q},x_1x_2^2$ give the type II
in \cite[example 3.2 (i)]{HK1}.
The following table lists all sets $M(k)$ and suitable 
values of $k$.
\begin{eqnarray*}
\begin{array}{ccc}
N&\{1\}&\emptyset\\ \hline
4q & 2q  & 1
\end{array}
\end{eqnarray*}
Therefore only the values of $L(k)=(-1)^{n-|M(k)|}\cdot \mu(k)$
and $\chi(k)$ for $k\in\{4q,2q,1\}$ are interesting.
\begin{eqnarray}
(L(4q),L(2q),L(1))=(2q+1,-(2q-1),1),\nonumber\\
(\chi(4q),\chi(2q),\chi(1))=(1,-1,1),\nonumber\\
D_{\bf w}=\Lambda_{4q}-\Lambda_{2q}+\Lambda_1.\label{5.2}
\end{eqnarray}
\end{examples}

K. Saito proposed the following conjecture.

\begin{conjecture}\label{t5.4}\cite[(3.13) and (4.2)]{Sa3}
Let $(w_1,...,w_n,1)$ be a normalized weight system
such that $\rho_{({\bf v},d)}\in\N_0[t]$
or (in general stronger) such that $(IS3)$ 
(from theorem \ref{t3.5}) holds for the reduced weight system.
Then $D_{\bf w}=\sum_{m\in\N}\nu(m)\cdot\Psi_m$ satisfies 
\begin{eqnarray}\label{5.3}
\nu(d_{\bf w})>0\textup{ or }\nu(\frac{d_{\bf w}}{2})>0,\\
\nu(d_{\bf w})>0\textup{ if all }w_j<\frac{1}{2},\label{5.4}
\end{eqnarray}
i.e. in the case $(IS3)$ the monodromy has eigenvalues
of order $d_{\bf w}$ or of order $\frac{d_{\bf w}}{2}$,
and if all $w_j<\frac{1}{2}$ it has eigenvalues of order
$d_{\bf w}$.
\end{conjecture}

Saito was not aware of the part of theorem \ref{t3.5}
saying that the condition $(C1)$ is sufficient for
$(IS3)$ (necessity is proved in \cite{Sa1}).
Probably therefore he gave in the conjecture in 
\cite[(3.13]{Sa3} the
characterization $\rho_{({\bf v},d)}\in\N_0[t]$,
which is in the cases $n\leq 3$ sufficient and
necessary for $(IS3)$ 
(\cite{AGV} \cite{Sa2} \cite[lemma 2.4]{HK1}).
In \cite[(4.2)]{Sa3} he gave the condition $(IS3)$.

He proved in \cite{Sa3} a result which implies the
conjecture for $n=3$. He also stated that it is true for $n=2$.

The following examples disprove the part \eqref{5.4}
of the conjecture for $n=4$. They can be extended easily
to $n\geq 5$.

\begin{examples}\label{t5.5}
Consider two curve singularities $D_{2^kq_1+1}$
and $D_{2^kq_2+1}$ with 
$k,q_1,q_2\in\N$ with $q_1$ and $q_2$ odd and
$\lcm(q_1,q_2)>\max(q_1,q_2)$.
Then their Thom-Sebastiani sum
$D_{2^kq_1+1}\otimes D_{2^kq_2+1}$ is a quasihomogeneous
singularity in $n=4$ variables with normalized weights
\begin{eqnarray*}
(\frac{1}{2^kq_1},\frac{2^kq_1-1}{2^{k+1}q_2},
\frac{1}{2^kq_2},\frac{2^kq_1-1}{2^{k+1}q_2})
\end{eqnarray*}
and $d_{\bf w}=2^{k+1}\lcm(q_1,q_2)$.
The divisor of the characteristic polynomial is
because of \eqref{5.2}
\begin{eqnarray*}
D_{\bf w}&=& (\Lambda_{2^{k+1}q_1}-\Lambda_{2^kq_1}+\Lambda_1)
\cdot (\Lambda_{2^{k+1}q_2}-\Lambda_{2^kq_2}+\Lambda_1)\\
&=& (2^{k+1}-2^k-2^k)\gcd(q_1,q_2)\Lambda_{2^{k+1}\lcm(q_1,q_2)}\\
&&+ 2^k\gcd(q_1,q_2)\Lambda_{2^k\lcm(q_1,q_2)}
+ \Lambda_{2^{k+1}q_1} +\Lambda_{2^{k+1}q_2}\\
&&-\Lambda_{2^kq_1} -\Lambda_{2^kq_2}+\Lambda_1.\\
&=& 2^k\gcd(q_1,q_2)\Lambda_{2^k\lcm(q_1,q_2)}
+\Lambda_{2^{k+1}q_1} +\Lambda_{2^{k+1}q_2}\\
&&-\Lambda_{2^kq_1} -\Lambda_{2^kq_2}+\Lambda_1.
\end{eqnarray*}
Part \eqref{5.4} of conjecture \ref{t5.3}
does not hold here.
\end{examples}

\begin{remarks}\label{t5.6}
In the examples \ref{t5.5} the part \eqref{5.3}
of the conjecture does hold. 
That part of the conjecture is still open.

We checked the tables of weight systems
of quasihomogeneous singularities in $n=4$
variables in \cite{HK2} up to $\mu=500$ for all
weight systems for which \eqref{5.4} does not hold.
There are 25 cases, and they are precisely
those Thom-Sebastiani sums 
$D_{2^kq_1+1}\otimes D_{2^kq_2+1}$ in the examples \ref{t5.4}
which satisfy $\mu\leq 500$. In 23 cases $k=1$, in 2 cases
$k=2$.

This indicates that for $n=4$ their might be no 
counterexamples to \eqref{5.3} and only the counterexamples
in example \ref{t5.5} to \eqref{5.4}.

{\bf Open problem 5:}
\begin{list}{}{}
\item[(a)] Prove or disprove the part \eqref{5.3} of conjecture
\ref{t5.4}. 
\item[(b)] Settle whether in the case $n=4$ the only
counterexamples to \eqref{5.4} are those in example
\ref{t5.5}.
\end{list}
\end{remarks}

\section{A conjecture on the orders of the
eigenvalues of the monodromy of a quasihomogeneous singularity}\label{s6}
\setcounter{equation}{0}

\noindent
Recall the definition \ref{t1.1} of the Orlik block
$(H_M,h_M)$ and of the group $\Aut_{S^1}(H_M,h_M)$ 
for a finite nonempty set $M\subset\N$.
The main result in \cite{He2} characterizes those
sets $M$ for which $\Aut_{S^1}(H_M,h_M)$ is as small
as possible in terms of conditions on the set $M$.
It is recalled below in theorem \ref{t6.2}.
The following definitions are needed.

\begin{definition}\label{t6.1}
Let $M\subset \Z_{\geq 1}$ be a finite set of positive integers.

\medskip
(a) A graph $\GG(M)=(M,E(M))$ is associated to it as follows. 
$M$ itself is the set of vertices. The edges in $E(M)$ are directed.
The set $E(m)$ is defined as follows. From a vertex $m_1\in M$ to
a vertex $m_2\in M$ there is no edge if at least one of the following 
two conditions holds:
\begin{list}{}{}
\item[(i)] 
$m_1/m_2$ is not a power of a prime number.
\item[(ii)]
An $m_3\in M-\{m_1,m_2\}$ with $m_2|m_3|m_1$ exists.
\end{list}
If $m_1/m_2$ is a power $p^k$ with $k\in\Z_{\geq 1}$ of a prime number $p$
and if no $m_3\in M-\{m_1,m_2\}$ with $m_2|m_3|m_1$ exists, 
then there is a directed edge from $m_1$ to $m_2$, 
which is additionally labelled with $p$. It is called a $p$-edge. 
Together such edges form the set $E(M)$ of all edges.

\medskip
(b) For any prime number $p$ the components of the graph 
$(M,E(M)-\{p\textup{-edges}\})$ which is obtained by deleting all
$p$-edges, are called the $p$-planes of the graph.
A $p$-plane is called a highest $p$-plane if no $p$-edge ends 
at a vertex
of the $p$-plane. A $p$-edge from $m_1$ to $m_2$ is called a
highest $p$-edge if no $p$-edge ends at $m_1$. 

\medskip
(c) A property $(T_p)$ for a prime number $p$ and a property
$(S_2)$ for the prime number 2:
\begin{eqnarray}\label{6.3}
(T_p) &:& \textup{The graph }\GG(M)
\textup{ has only one highest }p\textup{-plane.}\\
(S_2) &:& \textup{The graph }
(M,E(M)-\{\textup{highest }2\textup{-edges}\})
\nonumber \\
&& \textup{ has only 1 or 2 components.}\label{6.4}
\end{eqnarray}

(d) The least common multiple of the numbers
in $M$ is denoted $\lcm(M)\in\Z_{\geq 1}$. 
For any prime number $p$ denote 
\begin{eqnarray*}
l(m,p)&:=&\max(l\in\Z_{\geq 0}\, |\, p^l\textup{ divides }m)
\quad\textup{for any  }m\in\Z_{\geq 1},\\
l(M,p)&:=& \max(l(m,p)\, |\, m\in M)=l(\lcm(M),p).
\end{eqnarray*}
Then $m=\prod_{p\textup{ prime number}} p^{l(m,p)}.$
\end{definition}

\begin{theorem}\label{t6.2}\cite[theorem 1.2]{He2}
Let $M\subset\Z_{\geq 1}$ be a finite set of positive integers,
and let $(H_M,h_M)$ be its Orlik block. Then 
\begin{eqnarray}\label{6.6}
\Aut_{S^1}(H_M,h_M)=\{\pm h_M^k\, |\, k\in\Z\}
\end{eqnarray}
holds if and only if the graph $\GG(M)$ satisfies one of the following
two properties.
\begin{list}{}{}
\item[(I)] $\GG(M)$ is connected. It satisfies $(S_2)$.
It satisfies $(T_p)$ for any prime number $p\geq 3$.
\item[(II)] $\GG(M)$ has two components $M_1$ and $M_2$.
The graphs $\GG(M_1)$ and $\GG(M_2)$ are $2$-planes of $\GG(M)$
and satisfy $(T_p)$ for any prime number $p\geq 3$.
Furthermore
\begin{eqnarray}\label{6.7}
\gcd(\lcm(M_1),\lcm(M_2))&\in& \{1;2\},\\
l(M_1,2)&>&l(M_2,2)\in\{0;1\}.\label{6.8}
\end{eqnarray}
\end{list}
\end{theorem}

Motivated by Orlik's conjecture \ref{t1.3}, 
theorem \ref{t6.2}, and a search in the lists of
weight systems and associated divisors $D_{\bf w}$
in \cite{HK2}, here we propose the following conjecture.

\begin{conjecture}\label{t6.3} (= Conjecture \ref{t1.4})
For any quasihomogeneous singularity, each of the sets 
$M_1,...,M_{\nu_{max}}$ satisfies condition (I) in theorem \ref{t6.2}.
\end{conjecture}

\begin{remarks}\label{t6.4}
(i) {\bf Open problem 6:}
Prove conjecture \ref{t6.3} combinatorially (or disprove it by a counterexample).

\medskip
(ii) The conjecture is hard to deal with, because it 
requires to split the characteristic polynomial
into its elementary divisors (as also Orlik's conjecture).
It is not easy to extract from the formula for
$D_{\bf w}$, which is by the result of Milnor and Orlik
the divisor of the characteristic polynomial,
information about these elementary divisors.
This formula is rather nice in terms of the $\Lambda_m$
(though as a product, not a sum),
but the elementary divisors require to consider the $\Psi_m$.

\medskip
(iii) The example \ref{t6.5} (i) shows that the conditions
(I) and (II) together in theorem \ref{t6.2} do not behave 
well under tensor product. The example \ref{t6.5} (ii) shows
that condition (I) alone does not behave well under 
tensor product. This leads to the open problem 7.
It generalizes conjecture \ref{t6.3}.
A solution of problem 7 (a)+(b) would imply a positive solution of 
problem 6.

\medskip
(iv) {\bf Open problem 7:}
\begin{list}{}{}
\item[(a)] Find a natural condition for products $f$ of cyclotomic polynomials
which implies for any elementary divisor of $f$ condition
(I) in theorem \ref{t6.2}, and which is stable under
tensor product. 
\item[(b)] Prove that the characteristic polynomial of any quasihomogeneous
singularity satisfies this condition.
\end{list} 

\medskip
(iv) It would be desirable to have other ways to express
condition (I) in theorem \ref{t6.2}, e.g. 
in terms of the $\chi(m)$ of the divisor of a 
characteristic polynomial. But it is not clear
how they could look like.

\medskip
(v) Below conjecture \ref{t6.3} is therefore proved only
in a few cases, in theorem \ref{t6.9}.
The proofs use lemma \ref{t4.1} and corollary \ref{t4.3}.
\end{remarks}

\begin{examples}\label{t6.5}
(i) Consider $f_1:=\Phi_{12}\Phi_6^2\Phi_4^2\Phi_2$ and 
$f_2:=\Phi_5\Phi_1$. Then
\begin{eqnarray*}
f_1\otimes f_2 = \Phi_{60}\Phi_{30}^2\Phi_{20}^2\Phi_{12}\Phi_{10}
\Phi_6^2\Phi_4^2\Phi_2
\end{eqnarray*}
by \eqref{2.19}--\eqref{2.20}. 
Denote by $f_1=g_{1,1}\cdot g_{1,2}$ and $f_2=g_2$ and 
$f_1\otimes f_2=g_{\otimes,1}\cdot g_{\otimes,2}$ the decompositions into
elementary divisors and by $M_{1,1},M_{1,2},M_2,M_{\otimes,1}M_{\otimes,2}
\subset\N$ the corresponding sets. Then
\begin{tabular}{ll}
$M_{1,1}=\{12,6,4,2\}$ & satisfies condition (I), \\
$M_{1,2}=\{6,4\}$ & satisfies condition (II), \\
$M_{2}=\{5,1\}$ & satisfies condition (I), \\
$M_{\otimes,1}=\{60,30,20,12,10,6,4,2\}$ & satisfies condition (I), \\
$M_{\otimes,2}=\{30,20,6,4\}$ & satisfies neither (I) nor (II).
\end{tabular}

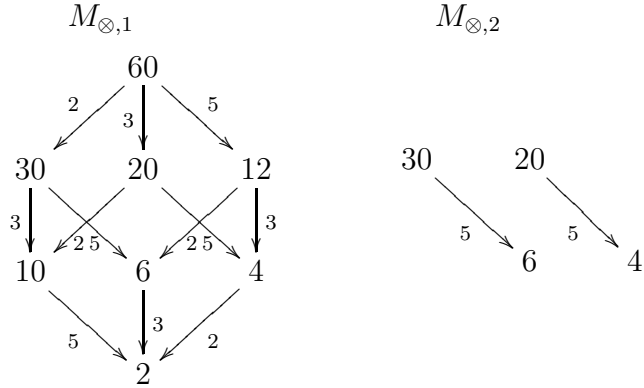
\begin{eqnarray*}
M_{\otimes,1} \hspace*{4cm}
M_{\otimes,2}  \\
\begin{xy} \hspace*{-3cm}
\xymatrix{ & 60 \ar[dl]_2 \ar[d]_3 \ar[dr]^5 & \\ 
30 \ar[d]_3 \ar[dr]_/1mm/2 & 20 \ar[dl]^/1mm/5 \ar[dr]_/1mm/2 & 
12 \ar[dl]^/1mm/5 \ar[d]^3 \\ 
10 \ar[dr]_{5} & 6 \ar[d]^3 & 4 \ar[dl]^{2}\\ & 2 &
}\hspace*{5cm} 
\xymatrix{ & & \\
 30 \ar[dr]_/1mm/5 & 20 \ar[dr]_/1mm/5 & \\  & 6 & 4
}  \end{xy}
\end{eqnarray*}

\medskip
(ii) Consider $f_1:=\Phi_7^2\Phi_3\Phi_1$ and 
$f_2:=\Phi_5^2\Phi_3\Phi_1$. Then
\begin{eqnarray*}
f_1\otimes f_2 = \Phi_{35}^4\Phi_{21}^2\Phi_{15}^2\Phi_7^2\Phi_5^2
\Phi_3^3\Phi_1^3
\end{eqnarray*}
by \eqref{2.19}--\eqref{2.21}. 
Denote by $f_1=g_{1,1}\cdot g_{1,2}$ and $f_2=g_{2,1}\cdot g_{2,2}$ and 
$f_1\otimes f_2=g_{\otimes,1}\cdot g_{\otimes,2}\cdot g_{\otimes,3}
\cdot g_{\otimes,4}$ the decompositions into
elementary divisors and by $M_{i,j}$ and $M_{\otimes,j}\subset\N$ 
the corresponding sets. Then
\begin{tabular}{ll}
$M_{1,1}=\{7,3,1\}$, & $M_{1,2}=\{7\}$, \\
$M_{2,1}=\{5,3,1\}$, & $M_{2,2}=\{5\}$, \\
$M_{\otimes,1}=M_{\otimes,2}=\{35,21,15,7,5,3,1\}$ & \textup{and }
$M_{\otimes,4}=\{35\}$ \\ & satisfy all condition (I), but\\
$M_{\otimes,3}=\{35,3,1\}$ & satisfies neither (I) nor (II).
\end{tabular}

\begin{eqnarray*}
M_{\otimes,1}=M_{\otimes,2} \hspace*{4cm}
M_{\otimes,3} \hspace*{2cm} M_{\otimes,4}\\
\begin{xy} \hspace*{-7cm}
\xymatrix{ 
21 \ar[d]_3 \ar[dr]_/1mm/5 & 35 \ar[dl]^/1mm/7 \ar[dr]_/1mm/5 & 
15 \ar[dl]^/1mm/7 \ar[d]^3 \\ 
7 \ar[dr]_{7} & 3 \ar[d]^3 & 5 \ar[dl]^{5}\\ & 1 &
}\hspace*{7cm} 
\xymatrix{ 
 35   \\ 3 \ar[d]^3 \\ 1} 
\hspace*{3cm} 
\xymatrix{ 
 35   \\  \\ } \end{xy}
\end{eqnarray*}
\end{examples}

\begin{remark}\label{t6.6}
(i) Lemma 8.2 in \cite{He1} gives the sufficient condition  in part (ii) 
for $\Aut_{S^1}(H_M,h_M)= \{\pm h_M^k\, |\, k\in\Z\}$.
It is a special case of condition (I) in theorem \ref{t6.2}.
It holds for many elementary divisors of characteristic polynomials
of quasihomogeneous singularities.
But the examples \ref{t6.7} (i)--(iii) give quasihomogeneous 
singularities where it does not hold for all elementary divisors
of the characteristic polynomial.

\medskip
(ii) A special case of condition (I) \cite[Lemma 8.2]{He1}:
$M$ contains a largest number $m_1$ such that $\GG(M)$
is a directed graph with root $m_1$. This implies $(T_p)$
for any $p$. Additionally, a chain of 2-edges exists which
connects all 2-planes. This implies $(S_2)$.
\end{remark}

\begin{examples}\label{t6.7}
(i) The weight system $({\bf w},1)=(\frac{1}{6},\frac{1}{10},
\frac{1}{15},1)$ satisfies $(C1)$ and $(IS3)$.
It is of type I (={\it Fermat type}) in the notation of
\cite[Example 3.2 (ii)]{HK1}. Here
\begin{eqnarray*}
D_{\bf w}&=& (\Lambda_6-\Lambda_1)(\Lambda_{10}-\Lambda_1)
(\Lambda_{15}-\Lambda_1)\\
&=&(2\Lambda_{30}-\Lambda_6-\Lambda_{10}+\Lambda_1)
(\Lambda_{15}-\Lambda_1)\\
&=& 20\Lambda_{30}+\Lambda_6+\Lambda_{10}+\Lambda_{15}-\Lambda_1
=\sum_{j=1}^{22}\divis g_j,
\end{eqnarray*}
with the elementary divisors $g_j$ with
\begin{eqnarray*}
\divis g_j &=& \Lambda_{30}\quad\textup{for }1\leq j\leq 20,\\
\divis g_{21} &=& \Lambda_6+\Lambda_{10}+\Lambda_{15}
-\Lambda_2-\Lambda_3-\Lambda_5+\Lambda_1\\
&=& \Psi_6+\Psi_{10}+\Psi_{15}+\Psi_2+\Psi_3+\Psi_5+\Psi_1,\\
\divis g_{22} &=& \Lambda_2+\Lambda_3+\Lambda_5-2\Lambda_1
=\Psi_2+\Psi_3+\Psi_5+\Psi_1.
\end{eqnarray*}
The sets $M_{21}=\{6,10,15,2,3,5,1\}$ for $g_{21}$ and
$M_{22}=\{2,3,5,1\}$ for $g_{22}$ satisfy condition (I) 
in theorem \ref{t6.2}, but not the stronger conditions 
in remark \ref{t6.6} (ii).

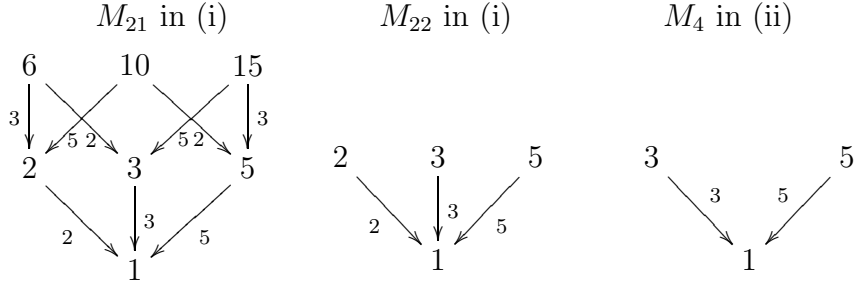
\begin{eqnarray*}
M_{21}\textup{ in (i)} \hspace*{2cm}
M_{22}\textup{ in (i)} \hspace*{2cm}
M_{4}\textup{ in (ii)} \\
\begin{xy} \hspace*{-7cm}
\xymatrix{
6 \ar[d]_3 \ar[dr]_/1mm/5 & 10 \ar[dl]^/1mm/2 \ar[dr]_/1mm/5 & 
15 \ar[dl]^/1mm/2 \ar[d]^3 \\ 
2 \ar[dr]_{2} & 3 \ar[d]^3 & 5 \ar[dl]^{5}\\ & 1 &
}\hspace*{4cm}
\xymatrix{ & & \\
2 \ar[dr]_{2} & 3 \ar[d]^3 & 5 \ar[dl]^{5}\\ & 1 &
}\hspace*{4cm}
\xymatrix{ & & \\
3 \ar[dr]^{3} & & 5 \ar[dl]_{5}\\ & 1 &
}\end{xy}
\end{eqnarray*}

\medskip
(ii) The weight system $({\bf w},1)=(\frac{2}{15},\frac{1}{5},
\frac{1}{3},1)$ satisfies $(C1)$ and $(IS3)$.
It is of type II in the notation of
\cite[Example 3.2 (ii)]{HK1}. Here
\begin{eqnarray*}
D_{\bf w}&=& (\frac{1}{2}\Lambda_{15}-\Lambda_1)
(\Lambda_5-\Lambda_1)(\Lambda_3-\Lambda_1)
=(2\Lambda_{15}-\Lambda_5+\Lambda_1)(\Lambda_3-\Lambda_1)\\
&=& 3\Lambda_{15}+\Lambda_5+\Lambda_3-\Lambda_1\\
&=& 3\Psi_{15}+4\Psi_5+4\Psi_3+4\Psi_1
=\sum_{j=1}^4\divis g_j,
\end{eqnarray*}
with the elementary divisors $g_j$ with
\begin{eqnarray*}
\divis g_j &=& \Lambda_{15}\quad\textup{for }1\leq j\leq 3,\\
\divis g_4 &=& 
\Psi_5+\Psi_3+\Psi_1.
\end{eqnarray*}
The set $M_4=\{5,3,1\}$ for $g_4$ satisfies condition (I) 
in theorem \ref{t6.2}, but not the stronger conditions 
in remark \ref{t6.6} (ii).

\medskip
(iii) The first of the examples \ref{t5.4} is $D_7\otimes D_{11}$
with $(k,q_1,q_2)=(1,3,5)$ and 
$({\bf w},1)=(\frac{1}{6},\frac{5}{12},\frac{1}{10},
\frac{9}{20},1)$. It satisfies $(C1)$ and is of type IV
in the notation of \cite[Example 3.2 (iii)]{HK1}. Here 
\begin{eqnarray*}
D_{\bf w}&=& 
2\Lambda_{30}+\Lambda_{12}+\Lambda_{20}-\Lambda_6
-\Lambda_{10}+\Lambda_1 =\sum_{j=1}^3\divis g_j,
\end{eqnarray*}
with the elementary divisors $g_j$ with
$\divis g_j =\sum_{m\in M_j}\Psi_m$ and 
\begin{eqnarray*}
M_1 &=& \{30,20,15,12,10,6,5,4,3,2,1\},\\
M_2 &=& \{30,15,10,6,5,4,3,2,1\},\\
M_3 &=& \{1\}.
\end{eqnarray*}
The sets $M_1$ and $M_2$ satisfy condition (I) 
in theorem \ref{t6.2}, but not the stronger conditions 
in remark \ref{t6.6} (ii).

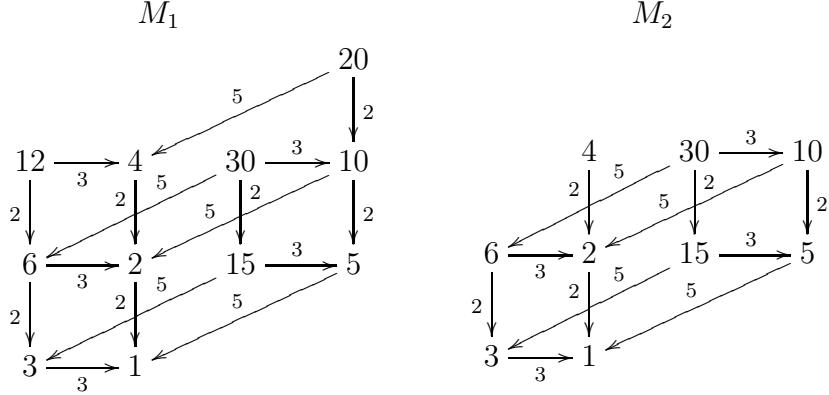
\begin{eqnarray*}
M_1 \hspace*{6cm} M_2\\
\begin{xy} \hspace*{-4cm}
\xymatrix{
 & & & 20 \ar[dll]_5 \ar[d]^2 \\
12 \ar[d]_2 \ar[r]_3 & 4 \ar[d]_/-2mm/2  & 
30 \ar[dll]_/-5mm/5 \ar[r]^3 \ar[d]^/-3mm/2 & 
10 \ar[dll]_/4mm/5 \ar[d]^2 \\
6 \ar[d]_2 \ar[r]_3 & 2 \ar[d]_/-2mm/2 & 
15 \ar[dll]_/-5mm/5 \ar[r]^3 & 5 \ar[dll]_5 \\
3 \ar[r]_3 & 1  & & \\
}\hspace*{6cm}
\xymatrix{ & & & \\
 & 4 \ar[d]_/-2mm/2 
& 30 \ar[dll]_/-5mm/5 \ar[r]^3 \ar[d]^/-3mm/2 & 
10 \ar[dll]_/4mm/5 \ar[d]^2 \\
6 \ar[d]_2 \ar[r]_3 & 2 \ar[d]_/-2mm/2 & 
15 \ar[dll]_/-5mm/5 \ar[r]^3 & 5 \ar[dll]_5 \\
3 \ar[r]_3 & 1  & & \\
}\end{xy}
\end{eqnarray*}
\end{examples}

\begin{lemma}\label{t6.8}
Suppose that numbers $k_1,...,k_l\in\N$ with
$k_{j}|k_{j-1}$ for $j=2,...,l$ are given.
Then the set $M\subset\N$ which is defined by 
$\Lambda_{k_1}-\Lambda_{k_2}+...
+(-1)^{l-1}\Lambda_{k_l}=\sum_{m\in M}\Psi_m$ is either empty or
satisfies the conditions in remark \ref{t6.6} (ii).
\end{lemma}

{\bf Proof:} 
We suppose that the set $M$ is not empty.
If $k_j=k_{j-1}$ for some $j\in\{2,3,...,l\}$, we can
drop $k_j$ and $k_{j-1}$. 
Therefore we can suppose $k_j<k_{j-1}$ for $j\in\{2,...,l\}$.
We have to prove the following two claims.

\medskip
{\bf Claim 1:}
{\it The graph $\GG(M)$ is a directed graph with root $k_1$.}

\medskip
{\bf Claim 2:}
{\it In $\GG(M)$ a chain of 2-edges exists which connects all
2-planes.}

\medskip
Proof of claim 1: The cases $l\in\{1,2\}$ are trivial.
Suppose $l\geq 3$. The proof uses induction over $l$.

Define the sets $M_1$ and $M_2$ by
$\Lambda_{k_1}-\Lambda_{k_2}=\sum_{m\in M_1}\Psi_m$
and $\sum_{j=3}^l(-1)^{j-1}\Lambda_{k_j}=\sum_{m\in M_2}\Psi_m$,
so that $M=M_1\dot\cup M_2$.
The graph $\GG(M_1)$ is obviously a directed graph with root
$k_1$. The graph $\GG(M_2)$ is by induction hypothesis a directed
graph with root $k_3$. For the proof of the claim it is sufficient
to show that the graph $\GG(M)$ contains a directed edge
from a vertex in $M_1$ to $k_3$. As $k_2<k_1$, a prime number
$q$ with $l(k_2,q)<l(k_1,q)$ exists. Then the number
\begin{eqnarray*}
m:= q^{l(k_1,q)}\cdot \prod_{p\ \textup{prime number},p\neq q}
p^{l(k_3,p)}
\end{eqnarray*}
is in $M_1$, and there is a directed edge from $m$
to $k_3$. \hfill $(\Box)$

\medskip
Useful for the proof of claim 2 will be {\bf Claim 3:}
{\it For any prime number $p$ and any $r\in\N_0$, the set
$M(p,r):= \{m\in M\, |\, l(m,p)=r\}$ is either empty or
a single $p-plane$. In the second case its graph
is a directed graph with a root.}

\medskip
Proof of claim 3: 
$M(p,r)=\{p^r\cdot m\, |\, m\in \www M(p,r)\}$ where
$\www M(p,r)$ is the support of the divisor
\begin{eqnarray*}
\sum_{j:\, l(k_j,p)\geq r}(-1)^{j-1}\Lambda_{\www{k_j}}
\quad\textup{with }\www{k_j}:=p^{-l(k_j,p)}\cdot k_j.
\end{eqnarray*}
If this divisor is not 0, claim 1 applies and gives claim 3.
\hfill $(\Box)$

\medskip
Proof of claim 2: Two cases will be distinguished.

{\bf 1st case}, for any odd $j\in\{1,...,l-1\}$ 
$\frac{k_j}{k_{j+1}}=2^{l(k_j,2)-l(k_{j+1},2)}:$ 
Then
\begin{eqnarray*}
M(2,r)&=& \left\{
\begin{array}{ll}
\emptyset & \textup{ if }r>l(k_1,2),\\
& \textup{ or if }l(k_j,2)\geq r>l(k_{j+1},2) \textup{ for an 
even }j,\\
& \textup{ or if }l(k_l,2)\geq r\textup{ and }
l\textup{ is even.}\\
\{2^r\cdot m &| \  m|\www {k_j}\}
\textup{ where }\www{k_j}:=2^{-l(k_j,2)}\cdot k_j\\
& \textup{ if }l(k_j,2)\geq r> l(k_{j+1},2)\textup{ for an odd }
j,\\
& \textup{ or if }l(k_l,2)\geq r\textup{ and }
j=l\textup{ is odd.}
\end{array}\right.
\end{eqnarray*}
Define
\begin{eqnarray*}
\www {k_{min}}:= \left\{\begin{array}{ll}
\www{k_l} & \textup{ if }l\textup{ is odd,}\\
\www{k_{l-1}} & \textup{ if }l\textup{ is even.}
\end{array}\right.
\end{eqnarray*}
Then the set 
$\{2^r\cdot \www{k_{min}}\, |\, M(2,r)\neq\emptyset\}$
is the set of vertices in $M$ of a chain of 2-edges which
connects all 2-planes.

\medskip
{\bf 2nd case}, a minimal odd $j\in\{1,2,...,l-1\}$
with $\frac{k_j}{k_{j+1}}\neq 2^{l(k_j,2)-l(k_{j+1},2)}$
exists: Then a prime number $p\geq 3$ with 
$l(k_j,p)> l(k_{j+1},p)$ exists.
And then $\{2^{r-l(k_j,2)}\cdot k_j\, |\, 0\leq r\leq l(k_j,2)\}
\subset M$. 
Therefore for $0\leq r\leq l(k_j,2)$ the set $M(2,r)$
is not empty.
For $r>l(k_j,2)$ the set $M(2,r)$ is as in the 1st case.
Thus the set $\{2^{r-l(k_j,2)}\cdot k_j\, |\, M(2,r)\neq 
\emptyset\}$ is the set of vertices in $M$ of a chain of
2-edges which connects all 2-planes.
\hfill$\Box$

\begin{theorem}\label{t6.9}
Conjecture \ref{t6.3} holds for the 
weight systems of quasihomogeneous singularities
of cycle type and of chain type.
It holds for all quasihomogeneous singularities 
in $n=2$ variables.
It holds for the quasihomogeneous singularities in $n=3$
variables which are of the types III, IV, V, VI and VII
in example 3.2 (ii) in \cite{HK1}
(see remark \ref{t6.10} for the types I and II).

In fact, in all these cases the set $M$ of each elementary
divisor satisfies even the stronger conditions  in 
remark \ref{t6.6} (ii).
\end{theorem}

{\bf Proof:} 
First consider the cycle type. Recall lemma \ref{t4.1},
and especially formula \eqref{4.5} for $D_{\bf w}$. 
It implies that all elementary divisors except one have
the divisor $\Lambda_{d/\gamma}$, and the last one
has the divisor $\Lambda_{d/\gamma}-\Lambda_1$ if
$n$ is odd, and it has the divisor $\Lambda_1$ if $n$ is even.
These divisors satisfy by lemma \ref{t6.8} the 
conditions in remark \ref{t6.6} (ii).

Next consider the chain type.
Recall corollary \ref{t4.3} and especially formula 
\eqref{4.17} for $D_{\bf w}$. 
It implies that any elementary divisor satisfies the
conditions in lemma \ref{t6.8}.
Therefore it satisfies the conditions in remark \ref{t6.6} (ii).

Now consider the case $n=2$. 
By example 3.2 (i) in \cite{HK1}, there are three types.
Type III is a cyle type. Type II is a chain type.
They are treated above.
Type I is the tensor product of two $A$-type singularities,
it is called Fermat type. In general, the tensor product is
difficult to deal with, but this case is fairly easy.
Here the weights are $(w_1,w_2)=(\frac{1}{t_1},\frac{1}{t_2})$,
and $D_{\bf w}$ is 
\begin{eqnarray*}
D_{\bf w}&=&(\Lambda_{t_1}-\Lambda_1)(\Lambda_{t_2}-\Lambda_1)\\
&=& \gcd(t_1,t_2)\Lambda_{\lcm(t_1,t_2)}-\Lambda_{t_1}
-\Lambda_{t_2}+\Lambda_1.
\end{eqnarray*}
The elementary divisors are as follows.
\begin{eqnarray*}
\textup{For }k\leq \gcd(t_1,t_2)-2:&&
\divis g_k=\Lambda_{\lcm(t_1,t_2)},\\
\textup{for }k= \gcd(t_1,t_2)-1:&&
\divis g_k=\Lambda_{\lcm(t_1,t_2)}
-\Lambda_{\gcd(t_1,t_2)}+\Lambda_1,\\
\textup{for }k= \gcd(t_1,t_2):&&
\divis g_k=\Lambda_{\lcm(t_1,t_2)}
-\Lambda_{t_1}-\Lambda_{t_2}+\Lambda_{\gcd(t_1,t_2)}.
\end{eqnarray*}
The divisors in the first two cases satisfy the conditions in 
lemma \ref{t6.8} and therefore the conditions in remark \ref{t6.6} (ii).

Consider the divisor $\divis g_k$ in the third case.
Suppose that $t_1\not| t_2$ and $t_2\not| t_1$, because else
$\divis g_k=0$. The set $M\subset\N$ with
$\divis g_k=\sum_{m\in M}\Psi_m$ 
is $$M=\{m\in\N\, |\, m|\lcm(t_1,t_2),m\not| t_1,m\not| t_2\}.$$
Obviously, its graph is a directed graph with root $\lcm(t_1,t_2)$.
This gives the first condition in remark \ref{t6.6} (ii).
For the second condition, we distinguish the following two cases.
Write $\www t_j=2^{-m(t_j,2)}\cdot t_j$, 
so that $t_j=2^{m(t_j,2)}\cdot \www t_j$.
Suppose $m(t_1,2)\geq m(t_2,2)$. Then $\www t_2\not|\www t_1$
and $\lcm(t_1,t_2)=2^{m(t_1,2)}\cdot\lcm(\www t_1,\www t_2)$.

1st case, $\www t_1\not| \www t_2$: Then the set
$\{2^r\cdot \lcm(\www t_1,\www t_2)\, |\, 0\leq r\leq m(t_1,2)\}$
is a subset of $M$ and is a chain of 2-edges which connects all 2-planes.

2nd case, $\www t_1|\www t_2$: Then $m(t_1,2)>m(t_2,2)$.
Then the set 
$\{2^r\cdot \www t_2\, |\, m(t_2,2)+1\leq r\leq m(t_1,2)\}$
is a subset of $M$ and is a chain of 2-edges which connects all 2-planes.

Now consider the case $n=3$. By example 3.2 (ii) in \cite{HK1},
there are seven types. Type V is a chain type, and type VII is a 
cycle type. They are treated above. The types III, IV and VI will
be treated in a similar way as the type I for $n=2$. 

Type III for $n=3$: weights 
${\bf w}=(\frac{1}{t_1},\frac{s_2}{t_2},\frac{s_3}{t_3})$
with 
\begin{eqnarray*}
w_j=\frac{1-w_1}{a_j},
\quad t_j=t_1\cdot\alpha_j\textup{ with }\alpha_j=\frac{a_j}{\gcd(a_j,t_1-1)}
\quad \textup{for }j=2,3,
\end{eqnarray*}
for some $a_2,a_3\in\N$. Write $\www\alpha:=\lcm(\alpha_2,\alpha_3)$.
Then
\begin{eqnarray*}
D_{\bf w}&=& (\frac{1}{s_1}\Lambda_{t_1}-\Lambda_1)
(\frac{1}{s_2}\Lambda_{t_2}-\Lambda_1)
(\frac{1}{s_3}\Lambda_{t_3}-\Lambda_1)\\
&=& (\Lambda_{t_1}-\Lambda_1)
\left(\frac{t_1\cdot\gcd(\alpha_2,\alpha_3)}{s_2\cdot s_3}\Lambda_{t_1\www\alpha}
-\frac{1}{s_2}\Lambda_{t_2}-\frac{1}{s_3}\Lambda_{t_3}+\Lambda_1\right)\\
&=& r_1\Lambda_{t_1\www\alpha}-r_2\Lambda_{t_2}-r_3\Lambda_{t_3}
+\Lambda_{t_1}-\Lambda_1\\
\textup{with} && r_1=\frac{t_1(t_1-1)\gcd(\alpha_2,\alpha_3)}{s_2s_3},\quad
r_2=\frac{t_1-1}{s_2},\quad r_3=\frac{t_1-1}{s_3}.
\end{eqnarray*}
Suppose (without loss of generality) that $r_2\leq r_3$.
The elementary divisors $g_k$ are as follows:
\begin{eqnarray*}
\textup{For }1\leq k\leq r_1-r_2-r_3: && \divis g_k =\Lambda_{t_1\www\alpha},\\
\textup{for }k=r_1-r_2-r_3+1: && \divis g_k = \Lambda_{t_1\www\alpha}
-\Lambda_{t_1\gcd(\alpha_2,\alpha_3)}\\
&&\hspace*{2cm} + \Lambda_{t_1} - \Lambda_1,\\
\textup{for }r_1-r_2-r_3+2\leq k\leq r_1-r_3: && 
\divis g_k = \Lambda_{t_1\www\alpha} - \Lambda_{t_1\gcd(\alpha_2,\alpha_3)},\\
\textup{for }r_1-r_3+1\leq k\leq r_1-r_2: && 
\divis g_k = \Lambda_{t_1\www\alpha} - \Lambda_{t_3},\\
\textup{for }r_1-r_2+1\leq k \leq r_1: && \divis g_k = \Lambda_{t_1\www\alpha}
-\Lambda_{t_2} - \Lambda_{t_3}\\
&& \hspace*{2cm} + \Lambda_{t_1\gcd(\alpha_2,\alpha_3)}.
\end{eqnarray*}
The divisors in the first four cases satisfy the conditions in lemma
\ref{t6.8} and therefore the conditions in remark \ref{t6.6} (ii).
The divisors in the fifth case is of the same type as the divisor
in the third case in type I for $n=2$. 

Type IV for $n=3$ is a sum of a 1 variable Fermat type and a 2 variable
cycle type. The weights are 
${\bf w}=(\frac{1}{t_1},\frac{s_2}{t_2},\frac{s_3}{t_3})$ with
\begin{eqnarray*}
\gamma=\gcd(a_2-1,a_2a_3-1)=\gcd(a_3-1,a_2a_3-1),\\ 
t_2=t_3=\frac{a_2a_3-1}{\gamma},\quad 
s_2=\frac{a_3-1}{\gamma},\quad 
s_3=\frac{a_2-1}{\gamma}
\end{eqnarray*}
for some $a_2,a_3\in\N_{\geq 2}$, by lemma \ref{t4.1}.
Write $\www t:=\gcd(t_1,t_2)$. 
Again by lemma \ref{t4.1}, $D_{\bf w}$ is the product
\begin{eqnarray*}
D_{\bf w} &=& (\Lambda_{t_1}-\Lambda_1)(\gamma\Lambda_{t_2}+\Lambda_1)\\
&=& \gamma\www t\Lambda_{\lcm(t_1,t_2)}-\gamma\Lambda_{t_2}+\Lambda_{t_1}
-\Lambda_1.
\end{eqnarray*}
The elementary divisors $g_k$ are as follows:
\begin{eqnarray*}
\textup{For }1\leq k\leq \gamma(\www t-1): && 
   \divis g_k =\Lambda_{\lcm(t_1,t_2)},\\
\textup{for }k=\gamma(\www t-1)+1: && 
   \divis g_k =\Lambda_{\lcm(t_1,t_2)} - \Lambda_{t_2}+\Lambda_{\www t}
   -\Lambda_1,\\
\textup{for }\gamma(\www t-1)+2\leq k\leq \gamma\www t: && 
   \divis g_k = \Lambda_{\lcm(t_1,t_2)} - \Lambda_{t_2},\\
\textup{for } k= \gamma\www t+1: && 
   \divis g_k = \Lambda_{t_1} - \Lambda_{\www t}.
\end{eqnarray*}
All these divisors satisfy the conditions in lemma \ref{t6.8} and 
therefore the conditions in remark \ref{t6.6} (ii).

Type VI for $n=3$ consists of a cycle such that one of its vertices
is the root of a 2 variable chain. The weights are 
${\bf w}=(\frac{s_1}{t_1},\frac{s_2}{t_2},\frac{s_3}{t_3})$ with
\begin{eqnarray*}
\gamma=\gcd(a_2-1,a_1a_2-1)=\gcd(a_1-1,a_1a_2-1),\ 
t_1=t_2=\frac{a_1a_2-1}{\gamma},\\ 
t_3=t_1\cdot\alpha\textup{ for some }\alpha\in\N,
\  s_1=\frac{a_2-1}{\gamma},\ 
s_2=\frac{a_1-1}{\gamma},
\end{eqnarray*}
for some $a_1,a_2\in \N_{\geq 2}$.
By lemma \ref{t4.1}, $D_{\bf w}$ is the product
\begin{eqnarray*}
D_{\bf w} &=& (\gamma\Lambda_{t_1}+\Lambda_1)(\frac{1}{s_3}\Lambda_{t_3}-\Lambda_1)\\
&=& r\Lambda_{t_3}-\gamma\Lambda_{t_1}-\Lambda_1
\quad\textup{with}\quad r=\frac{\gamma t_1+1}{s_3}.
\end{eqnarray*}
Observe $r\geq\gamma+1$, because $r-\gamma-1$ is the coefficient of 
$\langle 1\rangle$ in $D_{\bf w}$.
The elementary divisors $g_k$ are as follows:
\begin{eqnarray*}
\textup{For }1\leq k\leq r-\gamma-1: && 
   \divis g_k =\Lambda_{t_3},\\
\textup{for }k=r-\gamma: && 
   \divis g_k =\Lambda_{t_3} -\Lambda_1,\\
\textup{for }r-\gamma+1\leq k\leq r: && 
   \divis g_k = \Lambda_{t_3} - \Lambda_{t_1}.
\end{eqnarray*}
All these divisors satisfy the conditions in lemma \ref{t6.8} and 
therefore the conditions in remark \ref{t6.6} (ii).
\hfill$\Box$

\begin{remark}\label{t6.10}
In example 3.2 (ii) in \cite{HK1}, i.e. for $n=3$, type I is the Fermat type
with ${\bf w}=(\frac{1}{t_1},\frac{1}{t_2},\frac{1}{t_3}$),
and type II is the sum of a 1 variable Fermat type and
and 2 variable chain type, so 
${\bf w}=(\frac{1}{t_1},\frac{1}{t_2}, \frac{s_3}{t_3})$ with
$\frac{s_3}{t_3}=\frac{1-w_2}{a_3}$ and $a_3\geq 2$.
In both cases, a similar ansatz as in the proof of theorem  
\ref{t6.9} leads to an unpleasant multitude of different subcases.
The examples \ref{t6.7} (i)+(ii) show that in special cases of both types,
some elementary divisor does not satisfy the conditions in remark \ref{t6.6} (ii).
It does not seem worth to try to prove conjecture \ref{t6.3} in this way.
\end{remark}

\end{document}